\documentclass[11pt,a4paper]{amsart}

\theoremstyle{plain}
\newtheorem{thm}{Theorem}[section]
\newtheorem*{thm*}{Theorem}
\newtheorem{lem}[thm]{Lemma}
\newtheorem{prop}[thm]{Proposition}
\newtheorem{cor}[thm]{Corollary}
\theoremstyle{definition}

\newtheorem{exmp}[thm]{Example}
\newtheorem{conj}[thm]{Conjecture}

\theoremstyle{remark}
\newtheorem{rem}[thm]{Remark}

\usepackage{amsfonts, amsmath, amssymb, amsthm}
\usepackage{amsmath, amssymb}
\usepackage{bbm}

\usepackage{color}
\usepackage{graphicx}
\usepackage{url}
\usepackage{overpic}

\usepackage{psfrag}
\usepackage{paralist}
\usepackage[a4paper]{anysize}
\usepackage[all,ps]{xy}

\newcommand\Sing{\operatorname{sing}}
\newcommand{\Sym}{\operatorname{Sym}}

\newcommand{\s}{\sigma}
\newcommand{\g}{\gamma}
\newcommand{\Star}{\operatorname{st}}
\newcommand{\Link}{\operatorname{lk}}
\newcommand{\Cone}{\operatorname{cone}}
\newcommand{\Dim}{\operatorname{dim}}
\newcommand{\Codim}{\operatorname{codim}}

\newcommand{\mysim}{\hspace{-.14cm}\sim}
\newcommand{\B}{\mathbbm{D}}
\newcommand{\Id}{\operatorname{Id}}
\newcommand{\M}{\mathfrak{M}}
\newcommand{\m}{\mathfrak{m}}
\newcommand{\h}{\mathfrak{h}}

\newcommand\realPath[1]{|#1|}
\newcommand{\odd}{\operatorname{odd}}

\newcommand\CC{{\mathbb C}}
\newcommand\Sph{{\mathbb S}}
\newcommand\restr{\kern-1ex\mid}

\begin{document}

\title{Constructing Simplicial Branched Covers}
\author[Witte]{Nikolaus Witte}
 \address{Nikolaus Witte, Fachbereich Mathematik, AG~7, TU Darmstadt, 64289 Darmstadt, Germany}
 \email{witte@math.tu-berlin.de}
 \thanks{The authors is supported by Deutsche Forschungsgemeinschaft, DFG Research Group
   ``Polyhedral Surfaces''.}
 \date{\today}

\begin{abstract}
  Izmestiev and Joswig described how to obtain a simplicial covering space (the
  \emph{partial unfolding}) of a given simplicial complex, thus obtaining a simplicial branched
  cover~[Adv.\ Geom.\ 3(2):191-255, 2003].  We present a large class of branched covers which can be
  constructed via the partial unfolding.  In particular, for~$d\leq4$ every closed oriented PL
  $d$-manifold is the partial unfolding of some polytopal $d$-sphere.
\end{abstract}


\maketitle


\section{Introduction}
\label{sec:introduction}

\noindent
Branched covers are applied frequently in topology -- most prominently in the study, construction and
classification of closed oriented PL $d$-manifolds. 
First results are by
Alexander~\cite{alexander:NORS} in 1920, who observed that any closed
oriented PL $d$-manifold~$M$ is a branched cover of the $d$-sphere.
Unfortunately Alexander's proof does not allow for any (reasonable)
control over the number of sheets of the branched cover, nor over the
topology of the branching set: The number of sheets depends on the
size of some triangulation of~$M$ and the branching set is the
co-dimension 2-skeleton of the $d$-simplex.

However, in dimension $d\leq 4$, the situation is fairly well
understood.  By results of Hilden~\cite{hilden:TFBC} and
Montesinos~\cite{montesinos:MFBC} any closed oriented 3-manifold~$M$
arises as 3-fold simple branched cover of the 3-sphere branched over a
link.  In dimension four the situation becomes increasingly difficult.
First Piergallini~\cite{piergallini:FMS4} showed how to obtain any
closed oriented PL 4-manifold as a 4-fold branched cover of the
4-sphere branched over a transversally immersed
PL-surface~\cite{piergallini:FMS4}. Iori~\&
Piergallini~\cite{iori_piergallini:4MFNS} then improved the standing
result showing that the branching set may be realized locally flat if
one allows for a fifth sheet for the branched cover, thus proving a
long-standing conjecture by Montesinos~\cite{montesinos:CSR}. The
question as to whether any closed oriented PL 4-manifold can be obtained as
4-fold cover of the 4-sphere branched over a locally flat PL-surface
is still open.

For the partial unfolding and the construction of closed oriented
combinatorial 3-manifolds we recommend Izmestiev~\&
Joswig~\cite{izmestiev_joswig:BC}. Their
construction has recently been simplified significantly by Hilden,
Montesinos-Amilibia, Tejada~\& Toro~\cite{MR2218370}. For those able
to read German additional analysis and examples can be found
in~\cite{witte:DIPL}. The partial unfolding is implemented in the
software package \texttt{polymake}~\cite{polymake}.

This work has been greatly inspired by a paper of Hilden,
Montesinos-Amilibia, Tejada~\& Toro~\cite{MR2218370} and their bold
approach. However, the techniques developed in the following turn out
to differ substantially from the ideas in~\cite{MR2218370}, allowing
for stronger results in dimension three and generalization to
arbitrary dimensions.

\subsubsection*{Outline of the paper.}
After some basic definitions and notations the partial
unfolding~$\widehat{K}$ of a simplicial complex~$K$ is introduced. The
partial unfolding defines a projection $p:\widehat{K}\to K$ which is a
simplicial branched cover if~$K$ meets certain connectivity
assumptions. We define combinatorial models of key features of a
branched cover, namely the branching set and the monodromy
homomorphism.

Sections~\ref{sec:constructing_branched_covers}
and~\ref{sec:extending_triangulations} are related, yet self
contained.  The main result of this paper is presented in
Theorem~\ref{thm:constructing_branched_covers} and we give an explicit
construction of a combinatorial $d$-sphere~$S$, such that
$p:\widehat{S}\to S$ is equivalent to a given simple, $(d+1)$-fold
branched cover~$r:X\to\Sph^d$ (with some additional restriction for the
branching set of~$r$). Theorem~\ref{thm:constructing_branched_covers}
is then applied to the construction of closed oriented PL
$d$-manifolds as branched covers for~$d\leq 4$. The construction of~$S$ and
the proof of its correctness take up the entire
Section~\ref{sec:constructing_branched_covers}.

Finally, in Section~\ref{sec:extending_triangulations} we discuss how
to extend $k$-coloring of a subcomplex~$L\subset K$ of a simplicial
$d$-complex~$K$ to a $\max\{k,d+1\}$-coloring of a refinement~$K'$
of~$K$, such that~$L$ is again a subcomplex of~$K'$. Since~$K'$ is
constructed from~$K$ via finitely many stellar subdivisions of
edges, all properties invariant under these subdivisions are
preserved, e.g. polytopality, regularity, shellability, and others.
This improves an earlier result by Izmestiev~\cite{izmestiev:EOC}.

\subsection{Basic definitions and notations.}

A simplicial complex~$K$ is a \emph{combinatorial $d$-sphere} or
\emph{combinatorial $d$-ball} if it is piecewise linear homeomorphic
to the boundary of the $(d+1)$-simplex, respectively to the
$d$-simplex. Equivalently,~$K$ is a combinatorial $d$-sphere or
$d$-ball if there is a common refinement of~$K$ and the boundary of
the $(d+1)$-simplex, respectively the $d$-simplex. A simplicial
complex~$K$ is a \emph{combinatorial manifold} if the vertex link of
each vertex of~$K$ is a combinatorial sphere or a combinatorial ball.
A manifold~$M$ is PL if and only if~$M$ has a triangulation as a
combinatorial manifold. For an introduction to PL-topology see
Bj\"orner~\cite[Part~II]{bjoerner:TM}, Hudson~\cite{hudson:PLT}, and
Rourke~\& Sanderson~\cite{rourke_sanderson:PLT}.

A finite simplicial complex is \emph{pure} if all the inclusion
maximal faces, called the \emph{facets}, have the same dimension.  We
call a co-dimension 1-face of a pure simplicial complex~$K$a
\emph{ridge}, and the \emph{dual graph}~$\Gamma^*(K)$ of~$K$ has the
facets as its node set, and two nodes are adjacent if they share a
ridge. We denote the 1-skeleton of~$K$ by~$\Gamma (K)$, its
\emph{graph}.

Further it is often necessary to restrict ourselves to simplicial
complexes with certain connectivity properties: A pure simplicial
complex~$K$ is \emph{strongly connected} if its dual
graph~$\Gamma^*(K)$ is connected, and \emph{locally strongly
  connected} if the star $\Star_K(f)$ of~$f$ is strongly connected for
each face $f\in K$. If~$K$ is locally strongly connected, then
connected and strongly connected coincide. Further we call~$K$
\emph{locally strongly simply connected} if for each face $f\in K$
with co-dimension~$\geq2$ the link~$\Link_K(f)$ of~$f$ is simply
connected, and finally,~$K$ is \emph{nice} if it is locally strongly
connected and locally strongly simply connected.  Observe that
combinatorial manifolds are always nice.

Let $(\s_0,\s_1,\dots,\s_l)$ be an ordering of the facets of a pure
simplicial $d$-complex~$K$, and let~$D_i=\bigcup_{0\leq j\leq i}\s_j$
denote the union of the first~$i$ facets. We call the ordering
$(\s_0,\s_1,\dots,\s_l)$ a \emph{shelling} of~$K$ if $D_{i-1}\cap\s_i$
is a pure simplicial $(d-1)$-complex for $1\leq i\leq l$. If~$K$ is
the boundary complex of a simplicial $(d+1)$-polytope, then~$K$ admits
a shelling order which can be computed efficiently; see
Ziegler~\cite[Chapter~8]{ziegler:LOP}.

A simplicial complex obtained from a shellable complex by stellar
subdivision of a face is again shellable, a shellable sphere or ball
is a combinatorial sphere or ball, and for~$1\leq i\leq l$ the
intersection~$D_{i-1}\cap\s_i$ is a combinatorial $(d-1)$-ball (or
sphere). A shellable simplicial complex~$K$ is a wedge of balls or
spheres in general. If~$K$ is a manifold, then~$D_i$ is a
combinatorial $d$-ball (or sphere) for~$0\leq i\leq l$, and in
particular we have that~$D_{i-1}\cap\s_i$,~$D_i$, and hence~$K$ are
nice.  We call a face~$f\subset\s_i$ \emph{free} if $f\not\in
D_{i-1}$. In particular the (inclusion) minimal free faces describe
all free faces, and they are also called \emph{restriction sets} in
the theory of $h$-vectors of simplicial polytopes.

\subsection{The branched cover.}
\label{sec:bc}

The concept of a covering of a space~$Y$ by another space~$X$ is
generalized by Fox~\cite{fox:CSWS} to the notion of the branched
cover. Here a certain subset~$Y_{\Sing}\subset Y$ may violate the
conditions of a covering map. This allows for a wider application in
the construction of topological spaces.  It is essential for a
satisfactory theory of (branched) coverings to make certain
connectivity assumption for~$X$ and~$Y$. The spaces mostly considered
are Hausdorff, path connected, and locally path connected; see
Bredon~\cite[III.3.1]{bredon:TAG}.  Throughout we will restrict our
attention to coverings of manifolds hence they meet the connectivity
assumptions in~\cite{bredon:TAG}.

Consider a continuous map $h:Z\to Y$, and assume the
restriction~$h:Z\to h(Z)$ to be a covering. If~$h(Z)$ is dense in~$Y$
(and meets certain additional connectivity conditions) then there is a
surjective map $p:X\to Y$ with $Z\subset X$ and
$p\left|_{Z}\right.=h$. The map~$p$ is called a \emph{completion}
of~$h$, and any two completions $p:X\to Y$ and $p':X'\to Y$ are
equivalent in the sense that there exists a homeomorphism
$\varphi:X\to X'$ satisfying $p'\circ\varphi=p$ and
$\varphi\left|_Z\right.=\Id$.  The map~$p:X\to Y$ obtained this way is
a \emph{branched cover}, and we call the unique minimal subset
$Y_{\Sing}\subset Y$ such that the restriction of~$p$ to the preimage
of $Y\setminus Y_{\Sing}$ is a covering, the \emph{branching set}
of~$p$.  The restriction of~$p$ to $p^{-1}(Y\setminus Y_{\Sing})$ is
called the \emph{associated} covering of~$p$. If $h:Z\to Y$ is a
covering, then~$X=Z$, and~$p=h$ is a branched cover with empty
branching set.

\begin{exmp}\label{exmp:branched_cover}
  For $k\geq 2$ consider the map 
  \[
  p_k:\CC \to \CC: z \mapsto z^k.
  \]
  The restriction
  $p_k\left|_{\B^2}\right.$ is a $k$-fold branched cover $\B^2 \to
  \B^2$ with the single branch point~$\{0\}$.
\end{exmp}

We define the \emph{monodromy homomorphism}
\[
\m_p:\pi_1(Y\setminus Y_{\Sing},y_0)\to\Sym(p^{-1}(y_0))
\]
of a branched cover for a point $y_0\in Y\setminus Y_{\Sing}$ as the
monodromy homomorphism of the associated covering: If
$[\alpha]\in\pi_1(Y\setminus Y_{\Sing},y_0)$ is represented by a
closed path~$\alpha$ based at~$y_0$, then~$\m_p$ maps~$[\alpha]$ to
the permutation $(x_i\mapsto \alpha_i(1))$, where
$\{x_1,x_2,\dots,x_k\}=p^{-1}(y_0)$ is the preimage of~$y_0$ and
$\alpha_i:[0,1]\to X$ is the unique lifting of~$\alpha$ with
$p\circ\alpha_i=\alpha$ and $\alpha_i(0)=x_i$; see
Munkres~\cite[Lemma~79.1]{munkres:TOP} and Seifert~\&
Threlfall~\cite[\S~58]{seifert_threlfall:TOT}.  The \emph{monodromy
  group}~$\M_p$ is defined as the image of~$\m_p$.

Two branched covers $p:X\to Y$ and~$p':X'\to Y'$ are \emph{equivalent}
if there are homeomorphisms $\varphi:X\to X'$ and $\psi:Y\to Y'$ with
$\psi(Y_{\Sing})=Y'_{\Sing}$, such that $p'\circ\varphi=\psi\circ p$
holds.  The well known Theorem~\ref{thm:bc} is due to the uniqueness
of~$Y_{\Sing}$, and hence the uniqueness of the associated covering;
see Piergallini~\cite[p.~2]{piergallini:MBC}.

\begin{thm}\label{thm:bc}
  Let $p:X\to Y$ be a branched cover of a connected manifold~$Y$.
  Then~$p$ is uniquely determined up to equivalence by the branching
  set~$Y_{\Sing}$, and the monodromy homeomorphism~$\m_p$. In
  particular, the covering space~$X$ is determined up to homeomorphy.
\end{thm}

Let~$Y$ be a connected manifold and~$Y_{\Sing}$ a co-dimension~2
submanifold, possibly with a finite number of singularities. We call a
branched cover~$p$ \emph{simple} if the image~$\m_p(m)$ of any
meridial loop~$m$ around a non-singular point of the branching set is
a transposition in~$\M_p$. Note that the $k$-fold branched cover
$p_k\left|_{\B^2}\right. :\B^2\to\B^2$ presented in
Example~\ref{exmp:branched_cover} is not simple for $k\geq3$.

\subsection{The partial unfolding.}

The partial unfolding~$\widehat{K}$ of a simplicial complex~$K$ first
appeared in a paper by Izmestiev~\& Joswig~\cite{izmestiev_joswig:BC},
with some of the basic notions already developed
in Joswig~\cite{MR1900311}. The partial unfolding is closely related to the
complete unfolding, also defined in~\cite{izmestiev_joswig:BC}, but we
will not discuss the latter. The partial unfolding is a geometric
object defined entirely by the combinatorial structure of~$K$, and
comes along with a canonical \mbox{projection~$p:\widehat{K}\to K$.}

However, the partial unfolding~$\widehat{K}$ may not be a simplicial
complex. In general~$\widehat{K}$ is a pseudo-simplicial complex:
Let~$\Sigma$ be a collection of pairwise disjoint geometric simplices
with simplicial attaching maps for some pairs
$(\s,\tau)\in\Sigma\times\Sigma$, mapping a subcomplex of~$\s$
isomorphically to a subcomplex of~$\tau$. Identifying the subcomplexes
accordingly yields the quotient space~$\Sigma/\mysim$, which is called
a \emph{pseudo-simplicial complex} if the quotient map
$\Sigma\to\Sigma/\mysim$ restricted to any $\s\in\Sigma$ is bijective.
The last condition ensures that there are no self-identifications
within each simplex~$\s\in\Sigma$.

\subsubsection*{The group of projectivities.}
Let~$\s$ and~$\tau$ be neighboring facets of a finite, pure simplicial
complex~$K$, that is, $\s\cap\tau$ is a ridge. Then there is exactly
one vertex in~$\s$ which is not a vertex of~$\tau$ and vice versa,
hence a natural bijection~$\langle\s,\tau\rangle$ between the vertex
sets of~$\s$ and~$\tau$ is given by
\[
  \langle\s,\tau\rangle:V(\s)\to V(\tau) :
  v\mapsto
  \begin{cases}
    v &\text{if} \quad v\in\s\cap\tau\\
    \tau\setminus\s &\text{if} \quad v=\s\setminus\tau.
  \end{cases}
\]
The bijection~$\langle\s,\tau\rangle$ is called the \emph{perspectivity}
from~$\s$ to~$\tau$.

\begin{figure}[t]
  \centering
  \begin{overpic}[width=.57\textwidth]{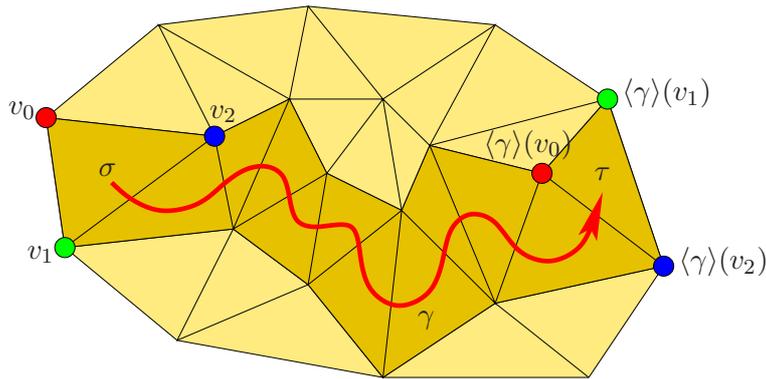}
    \put(13,32){$\mathbf{\s}$}
    \put(84.5,31.5){$\mathbf{\tau}$}
    \put(59,10.5){$\mathbf{\g}$}
    \put(0,41){$v_0$}
    \put(2.7,20){$v_1$}
    \put(29,40.5){$v_2$}
    \put(68.5,35.5){$\langle\g\rangle(v_0)$}
    \put(88.5,43){$\langle\g\rangle(v_1)$}
    \put(96.7,18.5){$\langle\g\rangle(v_2)$}
  \end{overpic}
  \caption{A projectivity from~$\s$ to~$\tau$ along the facet path~$\g$.\label{fig:projectivity}}
\end{figure}

A \emph{facet path} in~$K$ is a sequence $\g=(\s_0,\s_1,\dots,\s_k)$
of facets, such that the corresponding nodes in the dual
graph~$\Gamma^*(K)$ form a path, that is, $\s_i\cap\s_{i+1}$ is a
ridge for all $0\leq i<k$; see Figure~\ref{fig:projectivity}.  Now the
\emph{projectivity}~$\langle\g\rangle$ along~$\g$ is defined as the
composition of perspectivities $\langle\s_i,\s_{i+1}\rangle$,
thus~$\langle\g\rangle$ maps~$V(\s_0)$ to~$V(\s_k)$ bijectively via
\[
\langle\g\rangle = \langle\s_{k-1},\s_k\rangle\circ
\dots\circ\langle\s_1,\s_2\rangle\circ\langle\s_0,\s_1\rangle.
\]

We write $\g\,\delta=(\s_0,\s_1,\dots,\s_k,\s_{k+1},\dots,\s_{k+l})$
for the \emph{concatenation} of two facet paths $\g=(\s_0,
\s_1,\dots,\s_k)$ and $\delta=(\s_k,\s_{k+1},\dots,\s_{k+l})$, denote
by $\g^-=(\s_k,\s_{k-1},\dots,\s_0)$ the \emph{inverse path} of~$\g$,
and we call~$\g$ a \emph{closed} facet path based at~$\s_0$ if
$\s_0=\s_k$. The set of closed facet paths based at~$\s_0$ together
with the concatenation form a group, and a closed facet path~$\g$
based at~$\s_0$ acts on the set~$V(\s_0)$ via $\g\cdot
v=\langle\g\rangle(v)$ for $v\in V(\s_0)$. Via this action we obtain
the \emph{group of projectivities} $\Pi(K,\s_0)$ given by all
permutations $\langle\g\rangle$ of~$V(\s_0)$. The group of
projectivities is a subgroup of the symmetric group~$\Sym(V(\s_0))$ on
the vertices of~$\s_0$.

The projectivities along null-homotopic closed facet paths based
at~$\s_0$ generate the subgroup~$\Pi_0(K,\s_0)$ of~$\Pi(K,\s_0)$,
which is called the \emph{reduced group of projectivities}.  Finally,
if~$K$ is strongly connected then $\Pi(K,\s_0)$ and $\Pi(K,\s'_0)$,
respectively $\Pi_0(K,\s_0)$ and $\Pi_0(K,\s'_0)$, are isomorphic for
any two facets $\s_0,\s'_0\in K$. In this case we usually omit the
base facet in the notation of the (reduced) group of projectivities,
and write $\Pi(K)=\Pi(K,\s_0)$, respectively $\Pi_0(K)=\Pi_0(K,\s_0)$.

\subsubsection*{The odd subcomplex.}
Let~$K$ be locally strongly connected; in particular,~$K$ is pure. The
link of a co-dimension 2-face~$f$ is a graph which is connected
since~$K$ is locally strongly connected, and~$f$ is called \emph{even}
if the link $\Link_K(f)$ of~$f$ is 2-colorable (i.e. bipartite as a graph), and \emph{odd}
otherwise. We define the \emph{odd subcomplex} of~$K$ as all odd
co-dimension 2-faces (together with their proper faces), and denote it
by~$K_{\odd}$ (or sometimes~$\odd(K)$).

Assume that~$K$ is pure and admits a $(d+1)$-\emph{coloring} of its
graph~$\Gamma(K)$, that is, we assign one color of a set of~$d+1$
colors to each vertex of~$\Gamma(K)$ such that the two vertices of any
edge carry different colors. Observe that the $(d+1)$-coloring of~$K$
is minimal with respect to the number of colors, and is unique up to
renaming the colors if~$K$ is strongly connected. Simplicial complexes
that are $(d+1)$-colorable are called \emph{foldable}, since a
$(d+1)$-coloring defines a non-degenerated simplicial map of~$K$ to
the $(d+1)$-simplex. In other places in the literature foldable
simplicial complexes are sometimes called balanced.

\begin{lem}\label{lem:foldable_and_odd}
  The odd subcomplex of a foldable simplicial complex~$K$ is empty,
  and the group of projectivities $\Pi(K,\s_0)$ is trivial. In
  particular we have $\langle\gamma\rangle=\langle\delta\rangle$ for
  any two facet paths~$\gamma$ and~$\delta$ from~$\s$ to~$\tau$ for
  any two facets $\s,\tau\in K$.
\end{lem}

We leave the proof to the reader. As we will see in
Theorem~\ref{thm:IJ_1} the odd subcomplex is of interest in particular
for its relation to $\Pi_0(K,\s_0)$ of a nice simplicial complex~$K$.

Consider a geometric realization~$|K|$ of~$K$. To a given facet
path~$\g=(\s_0,\s_1,\dots,\s_k)$ in~$K$ we associate a (piecewise
linear) path~$\realPath{\g}$ in~$|K|$ by connecting the barycenter
of~$\s_i$ to the barycenters of $\s_i\cap\s_{i-1}$ and
$\s_i\cap\s_{i+1}$ by a straight line for $1\leq i<k$, and connecting
the barycenters of~$\s_0$ and $\s_0\cap\s_1$, respectively $\s_k$ and
$\s_k\cap\s_{k-1}$. A projectivity \emph{around} a co-dimension
2-face~$f$ is a projectivity along a facet path $\g\,\delta\,\g^-$,
where~$\delta$ is a closed facet path in $\Star_K(f)$ (based at some
facet $\s\in\Star_K(f)$) such that $\realPath{\g}$ is homotopy
equivalent to the boundary of a transversal disc
around~$|f|\subset|\Star_K(f)|$, and~$\g$ is a facet path from~$\s_0$
to~$\s$. The path $\g\,\delta\,\g^-$ is null-homotopic since~$K$ is
locally strongly simply connected.

\begin{thm}[{Izmestiev \& Joswig~\cite[Theorem~3.2.2]{izmestiev_joswig:BC}}]
  \label{thm:IJ_1}
  The reduced group of projectivities~$\Pi_0(K,\s_0)$ of a nice
  simplicial complex~$K$ is generated by projectivities around the odd
  co-dimension 2-faces. In particular,~$\Pi_0(K,\s_0)$ is generated by
  transpositions.
\end{thm}

The fundamental group~$\pi_1(|K|\setminus|K_{\odd}|,y_0)$ of a nice
simplicial complex~$K$ is generated by paths~$\realPath{\g}$,
where~$\g$ is a closed facet path based at~$\s_0$, and~$y_0$ is the
barycenter of~$\s_0$;
see~\cite[Proposition~A.2.1]{izmestiev_joswig:BC}.  Furthermore, due
to Theorem~\ref{thm:IJ_1} we have the group homomorphism
\[
\h_K:\pi_1(|K|\setminus|K_{\odd}|,y_0)\to\Pi(K,\s_0):[\realPath{\g}]\mapsto\langle\g\rangle,
\]
where $[\realPath{\g}]$ is the homotopy class of the
path~$\realPath{\g}$ corresponding to a facet path~$\g$.

\subsubsection*{The partial unfolding}
\begin{figure}[htbp]
  \centering
  \begin{overpic}[width=.65\textwidth]{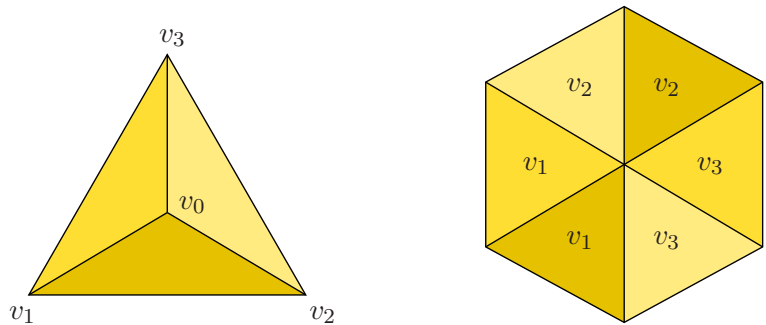}
    \put(22.5,18){$v_0$}
    \put(1,4){$v_1$}
    \put(39,4){$v_2$}
    \put(20,39.3){$v_3$}
    \put(66,23){$v_1$}
    \put(71.5,33){$v_2$}
    \put(88,23){$v_3$}
    \put(71.5,13.5){$v_1$}
    \put(82.5,33){$v_2$}
    \put(82.5,13.5){$v_3$}
  \end{overpic}
  \caption{The starred triangle and its partial unfolding. The complex on the
    right is the non-trivial connected component
    of the partial unfolding, indicated by the labeling of the facets
    by the vertices~$v_1$,~$v_2$, and~$v_3$. The second connected
    component is a copy of the starred triangle with all facets
    labeled~$v_0$; see also Example~\ref{exmp:branched_cover} for~$k=2$.
    \label{fig:part_unf}}
\end{figure}

Let~$K$ be a pure simplicial $d$-complex and set~$\Sigma$ as the set
of all pairs $(|\s|,v)$, where~$\s\in K$ is a facet and~$v\in\s$ is a
vertex.  Thus each pair $(|\s|,v)\in\Sigma$ is a copy of the geometric
simplex~$|\s|$ labeled by one of its vertices.  For neighboring
facets~$\s$ and~$\tau$ of~$K$ we define the equivalence
relation~$\sim$ by attaching $(|\s|,v)\in\Sigma$ and
$(|\tau|,w)\in\Sigma$ along their common ridge~$|\s\cap\tau|$ if
$\langle\s,\tau\rangle(v)=w$ holds. Now the \emph{partial
  unfolding}~$\widehat{K}$ is defined as the quotient space
$\widehat{K} = \Sigma/\mysim$. The projection $p:\widehat{K}\to K$ is
given by the factorization of the map $\Sigma\to K:(|\s|,v)\mapsto
\s$; see Figure~\ref{fig:part_unf}.

The partial unfolding of a strongly connected simplicial complex is
not strongly connected in general. We denote by
$\widehat{K}_{(|\s|,v)}$ the connected component containing the
labeled facet $(|\s|,v)$. Clearly,
$\widehat{K}_{(|\s|,v)}=\widehat{K}_{(|\tau|,w)}$ holds if and only if
there is a facet path~$\g$ from~$\s$ to~$\tau$ in~$K$ with
$\langle\g\rangle(v)=w$. It follows that the connected components
of~$\widehat{K}$ correspond to the orbits of the action
of~$\Pi(K,\s_0)$ on~$V(\s_0)$.  Note that each connected component of
the partial unfolding is strongly connected and locally strongly
connected~\cite[Satz~3.2.2]{witte:DIPL}. Therefore we do not
distinguish between connected and strongly connected components of the
partial unfolding.

The problem that the partial unfolding~$\widehat{K}$ may not be a
simplicial complex can be addressed in several ways. Izmestiev~\&
Joswig~\cite{izmestiev_joswig:BC} suggest barycentric subdivision
of~$\widehat{K}$, or anti-prismatic subdivision of~$K$.  A more
efficient solution (with respect to the size of the resulting
triangulations) is given in~\cite{witte:DIPL}.

\subsection{The partial unfolding as a branched cover.}

As preliminaries to this section we state two theorems by
Fox~\cite{fox:CSWS} and Izmestiev~\&
Joswig~\cite{izmestiev_joswig:BC}. Together they imply that under the
``usual connectivity assumptions'' the partial unfolding of a
simplicial complex is indeed a branched cover as suggested in the
heading of this section.

\begin{thm}[{Izmestiev \& Joswig~\cite[Theorem 3.3.2]{izmestiev_joswig:BC}}]
  Let~$K$ be a nice simplicial complex. Then the restriction of
  $p:\widehat{K}\to K$ to the preimage of the complement of the odd
  subcomplex is a simple covering.
\end{thm}

\begin{thm}[{Fox~\cite[p.~251]{fox:CSWS}; Izmestiev~\&
    Joswig~\cite[Proposition 4.1.2]{izmestiev_joswig:BC}}]
  Let ~$J$ and~$K$ be nice simplicial complexes and let $f:J\to K$ be
  a simplicial map. Then the map~$f$ is a simplicial branched cover if
  and only if
  \[
  \Codim K_{\Sing}\geq2.
  \]
\end{thm}

Since the partial unfolding of a nice simplicial complex is nice
Corollary~\ref{cor:unf_are_covers} follows immediately.

\begin{cor}\label{cor:unf_are_covers}
  Let~$K$ be a nice simplicial complex. The projection
  $p:\widehat{K}\to K$ is a simple branched cover with the odd
  subcomplex~$K_{\odd}$ as its branching set.
\end{cor}

For the rest of this section let~$K$ be a nice simplicial complex and
let~$y_0$ be the barycenter of a fixed facet~$\s_0\in K$. The
projection $p:\widehat{K}\to K$ is a branched cover with
$K_{\Sing}=K_{\odd}$ by Corollary~\ref{cor:unf_are_covers}, and
Izmestiev~\& Joswig~\cite{izmestiev_joswig:BC} proved that there is a
bijection $\imath:p^{-1}(y_0)\to V(\s_0)$ that induces a group
isomorphism $\imath_*:\Sym(p^{-1}(y_0))\to\Sym(V(\s_0))$ such that the
following Diagram~\eqref{equ:partial_unf} commutes.
\begin{equation}\label{equ:partial_unf}
\xymatrix{ 
  \pi_1(|K|\setminus|K_{\odd}|,y_0) \ar[dr]^{\h_K} \ar[d]^{\m_p} & \\
  \M_p \ar[r]_{{\imath_*}
  } & \Pi(K,\s_0) \\
}
\end{equation}

Let $r:X\to Y$ be a branched cover and assume that there is a
homomorphism of pairs $\varphi:(Y,Y_{\Sing})\to(|K|,|K_{\odd}|)$, that
is, $\varphi:Y\to|K|$ is a homomorphism with
$\varphi(Y_{\Sing})=|K_{\odd}|$.  Then Theorem~\ref{thm:partial_unf}
gives sufficient conditions for $p:\widehat{K}\to K$ and $r:X\to Y$ to
be equivalent branched covers.  It is the key tool in the proof of the
main Theorem~\ref{thm:constructing_branched_covers} in
Section~\ref{sec:constructing_branched_covers}.

\begin{thm}\label{thm:partial_unf}
  Let~$K$ be a nice simplicial complex and let $r:X\to Y$ be a
  (simple) branched cover. Further assume that there is a homomorphism
  of pairs $\varphi:(Y,Y_{\Sing})\to(|K|,|K_{\odd}|)$, and let $y_0\in
  Y$ be a point such that~$\varphi(y_0)$ is the barycenter of~$|\s_0|$
  for some facet $\s_0\in K$. The branched covers $p:\widehat{K}\to K$
  and $r:X\to Y$ are equivalent if there is a bijection
  $\iota:r^{-1}(y_0)\to V(\s_0)$ that induces a group isomorphism
  $\iota_*:\M_r\to\Pi(K,\s_0)$ such that the diagram
  \begin{equation}\label{equ:partial_unf_and_bc}
    \xymatrix{
      \pi_1(Y\setminus Y_{\Sing},y_0) \ar[d]^{\m_r} \ar[r]^{\varphi_*} &
      \pi_1(|K|\setminus|K_{\odd}|,\varphi(y_0)) \ar[d]^{\h_K}\\
      \M_r \ar[r]^{\iota_*} & \Pi(K,\s_0)
    }
  \end{equation}
  commutes. In particular, we have $\widehat{K}\cong X$.
\end{thm}
  
\begin{proof}
  Corollary~\ref{cor:unf_are_covers} ensures that $p:\widehat{K}\to K$
  is indeed a branched cover, and commutativity of
  Diagram~\eqref{equ:partial_unf} and
  Diagram~\eqref{equ:partial_unf_and_bc} proves commutativity of their
  composition:
  \[
  \xymatrix{
    \pi_1(Y\setminus Y_{\Sing},y_0) \ar[d]^{\m_r} \ar[r]^{\varphi_*} &
    \pi_1(|K|\setminus|K_{\odd}|,\varphi(y_0)) \ar[d]^{\h_K}
    \ar[dr]^{\m_p} &\\
    \M_r \ar[r]^{\iota_*} 
    & \Pi(K,\s_0) & \M_p \ar[l]_{\imath_*}
  }
  \]
  Theorem~\ref{thm:bc} completes the proof.
\end{proof}


\section{Constructing Branched Covers}
\label{sec:constructing_branched_covers}

\noindent
Throughout this section let $r:X\to\Sph^d$ be a branched cover of the
$d$-sphere with branching set~$F$.  The main objective is to 
give a large class of branched covers~$r$, such that there is a combinatorial sphere~$S$
with $p:\widehat{S}\to S$ equivalent to~$r$ as a branched cover. In
particular this implies the existence of a homeomorphism of pairs
$\varphi:(\Sph^d,F)\to(|S|,|S_{\odd}|)$.  Note that by the nature of
the partial unfolding and the projection~$p:\widehat{S}\to S$ any branched cover~$r$ equivalent to~$p$ has to be
simple and $(d+1)$-fold. A theorem similar to
Theorem~\ref{thm:constructing_branched_covers} may easily be
formulated for branched covers of $d$-balls.

\enlargethispage{\baselineskip}
Recall that we associate to a facet path~$\g$ in~$S$ the (realized)
path $\realPath{\g}$ in~$|S|$, and that the square brackets denote the
homotopy class of a closed path. Thus we write
$\m_r([\varphi^{-1}(\realPath{\g})])$ for the image of an element in
$\pi_1(\Sph^d\setminus F,y_0)$ represented by the closed
path~$\varphi^{-1}(\realPath{\g})$, which in turn is obtained from a
closed facet path~$\g$ based at some facet~$\s_0\in S$ with barycenter
$\varphi(y_0)$ by first considering its
realization~$\realPath{\gamma}$ and then its preimage under~$\varphi$.

\begin{thm}\label{thm:constructing_branched_covers}
  For $d\geq2$ let $r:X\to\Sph^d$ be a $(d+1)$-fold, simple branched
  cover of the $d$-sphere, and assume that the branching set~$F$
  of~$r$ can be embedded via a homeomorphism $\varphi:\Sph^d\to|S'|$
  into the co-dimension $2$-skeleton of a shellable simplicial
  $d$-sphere~$S'$.  Then there is a shellable simplicial $d$-sphere~$S$,
  such that $p:\widehat{S}\to S$ is a branched cover equivalent
  to~$r$. Further more, the $d$-sphere~$S$ can be obtained from~$S'$
  by a finite series of stellar subdivision of edges. If~$S'$ is the
  boundary of a simplicial $(d+1)$-polytope then also~$S$ is the
  boundary of a simplicial $(d+1)$-polytope.
\end{thm}

To make the proof of Theorem~\ref{thm:constructing_branched_covers}
more digestible we first give the (algorithmical) back-bone of the
proof
and defer some of the more technical aspects to the
Lemmas~\ref{lem:equivalence},~\ref{lem:non_free_edge},
and~\ref{lem:terminates}. Fix a point~$y_0\in\Sph^d\setminus F$ and
we may assume~$\varphi(y_0)$ to be the barycenter of some
facet~$\s_0\in S'$ and $|\s_0|\cap\varphi(F)=\emptyset$ to hold.
Further fix a bijection~$\imath$ between the
preimage~$\{x_0,x_1,\dots,x_d\}=r^{-1}(y_0)$ of~$y_0$ and the vertices
of~$\s_0$, and color the vertices of~$\s_0$ via~$\imath$ by the
elements in~$r^{-1}(y_0)$.

\begin{figure}[htbp]
  \centering
  \begin{overpic}[height=6.3cm]{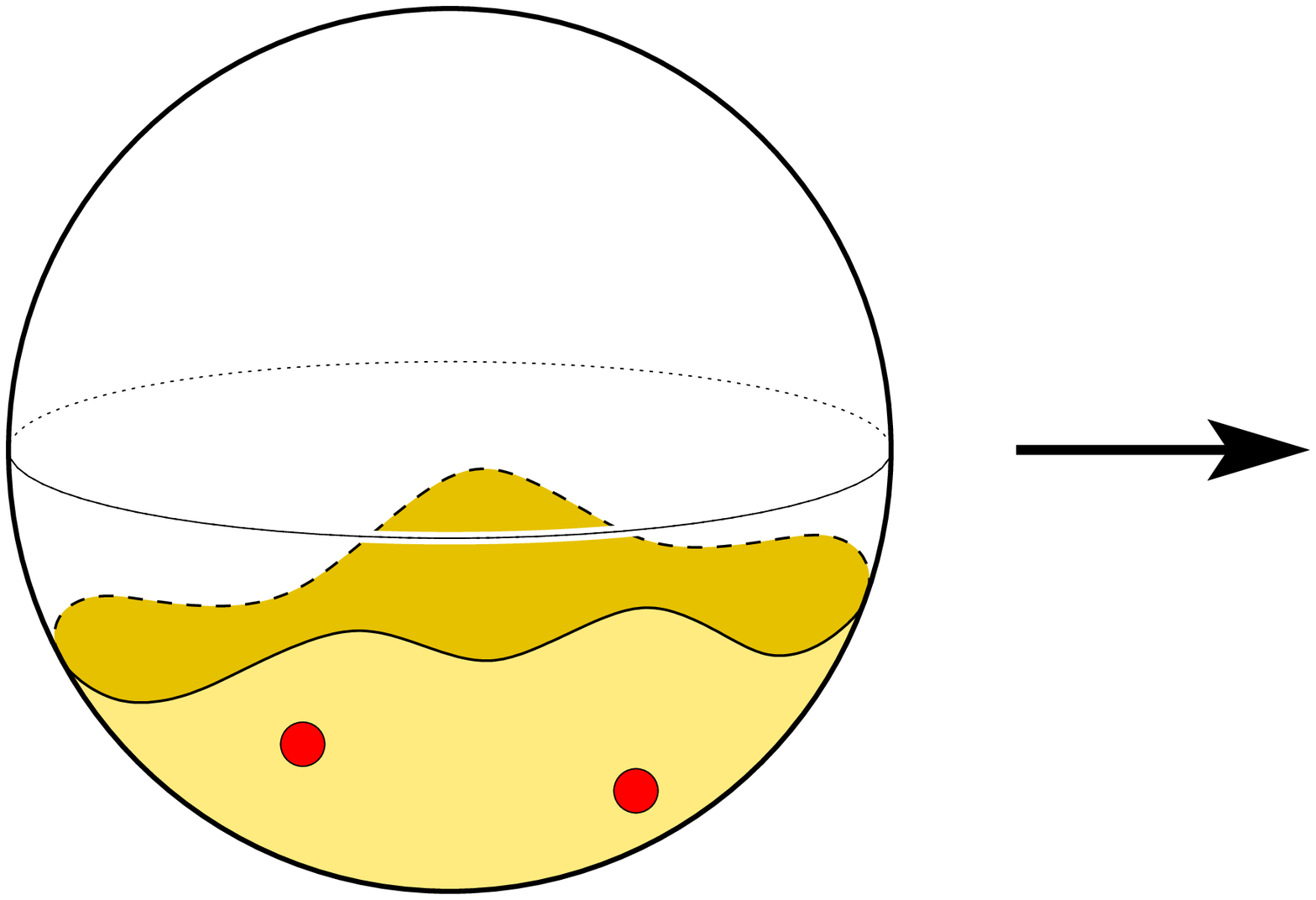}
    \put(83,38.5){\Large$\mathbf{\varphi}$}
  \end{overpic}
  \includegraphics[height=6.3cm]{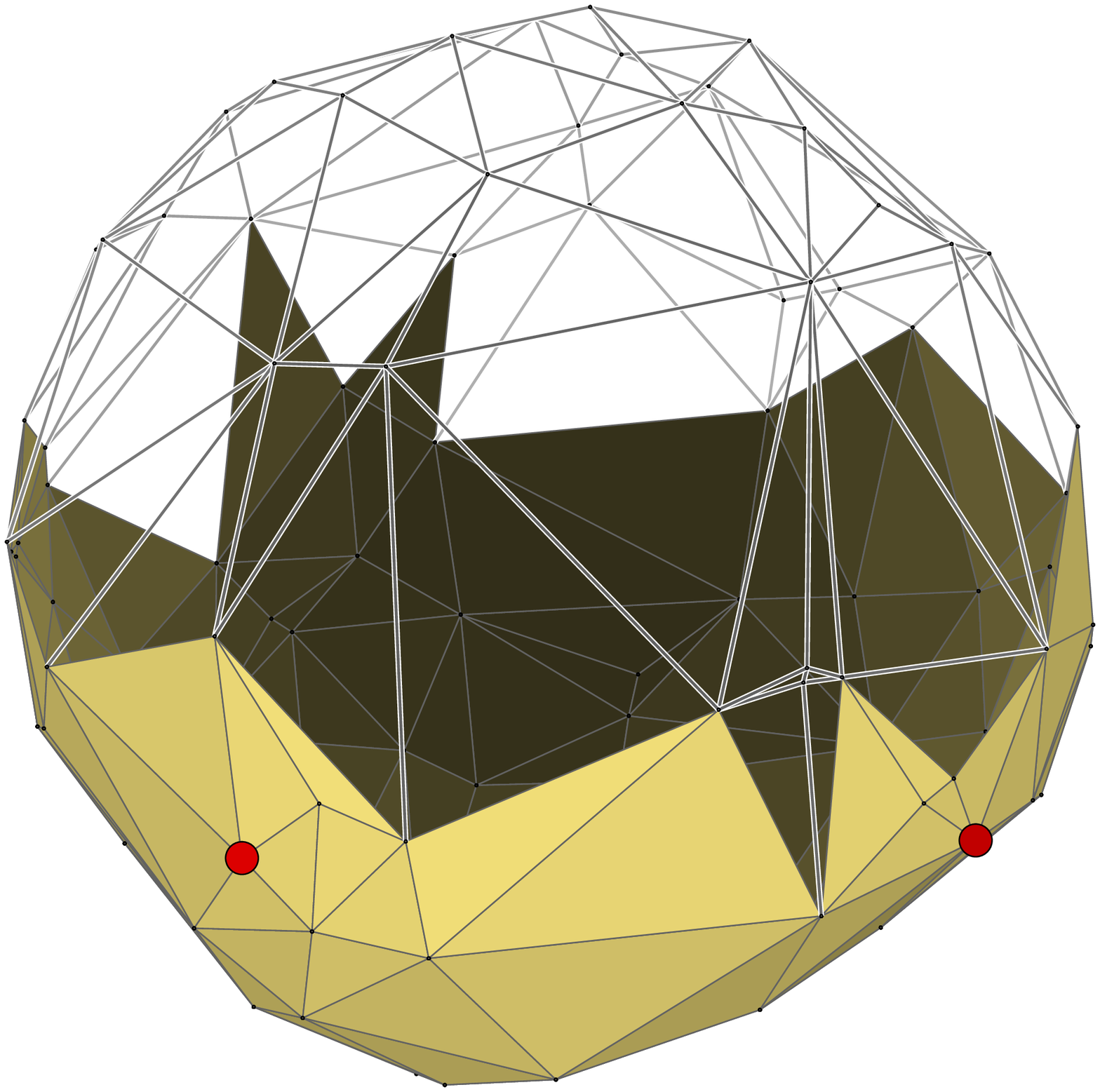}
  \caption{
    The base space of the branched cover~$r:X\to\Sph^2$ (left)
    and a 
    polytopal
    $2$-sphere~$S_i$ with marked beginning~$(\s_j)_{0\leq j\leq l_i}$ of a
    shelling (right). On the left the preimage of~$D_i=\bigcup_{0\leq j\leq l_i}\s_j$
    under the homomorphism~$\varphi:\Sph^2\to|S_i|$
    is shaded and the branching set is marked. The odd subcomplex of~$D_i$ is
    marked on the right. The branched covers~$r:X\to\Sph^2$
    (restricted to~$\varphi^{-1}(|D_i|)$) and~$\widehat{D_i}\to
    D_i$ are equivalent.
    \label{fig:shelling}}
\end{figure}

The $d$-sphere~$S$ is constructed in a finite series
$(S'=S_0,S_1,\dots,S_l=S)$ of shellable $d$-spheres, and each $d$-sphere~$S_i$
comes with a shelling of its facet with marked beginning
$(\s_{i,0},\s_{i,1},\dots,\s_{i,l_i})$.  The complex~$S_{i+1}$ is
obtained from~$S_i$ by (possibly) subdividing~$\s_{i,l_i+1}$ in a
finite series of stellar subdivisions of edges not contained in
any~$\s_{i,j}$ for $0\leq j\leq l_i$. Thus we may choose the shelling
of~$S_{i+1}$ such that it extends
$(\s_{i,0},\s_{i,1},\dots,\s_{i,l_i})$ and we
denote the marked beginning of the shelling of~$S_i$ simply by
$(\s_0,\s_1,\dots,\s_{l_i})$. 

Let $D_i=\bigcup_{0\leq j\leq l_i}\s_j$ then the
main idea of the proof of
Theorem~\ref{thm:constructing_branched_covers} is to construct~$S_i$
such that the branched covers $r:X\to\Sph^d$ (restricted
to~$\varphi^{-1}(|D_i|)$) and~$\widehat{D_i}\to D_i$ are
equivalent. To this end we prove that~$\varphi$ restricted to~$\varphi^{-1}(|D_i|)$ is a
homomorphism of pairs $(\varphi^{-1}(|D_i|),F\cap \varphi^{-1}(|D_i|))\to(|D_i|,|\odd(D_i)|)$
and that the
following Diagram~\eqref{eq:constructing_branched_covers} commutes; see Figure~\ref{fig:shelling}.
\begin{equation}\label{eq:constructing_branched_covers}
  \xymatrix{
    \pi_1(\varphi^{-1}(|D_i|)\setminus F,y_0) \ar[r]^{\varphi_*} \ar[d]^{\m_r} &
    \pi_1(|D_i|\setminus|\odd(D_i)|,\varphi(y_0)) \ar[d]^{\h_{D_i}}\\
    \M_r \ar[r]^{\imath_*} & \Pi(D_i,\s_0)
  }
\end{equation}
Commutativity of Diagram~\eqref{eq:constructing_branched_covers} is
obtained by ensuring that for each closed facet path~$\g$
in~$D_i$ (which is not a facet path in~$D_{i-1})$ the projectivity~$\langle\g\rangle$ acts
on~$V(\s_0)$ as~$m_r([\varphi^{-1}(\realPath{\g})])$ acts on~$r^{-1}(y_0)$.

The pair $(S_{i+1},(\s_j)_{0\leq j\leq l_{i+1}})$ is constructed from
the pair $(S_i,(\s_j)_{0\leq j\leq l_i})$ as follows.
Let~$\s=\s_{l_i+1}$ be the first facet in the shelling of~$S_i$ not
contained in $D_i$, let~$\g$ be a facet path in~$D_i\cup\s$
from~$\s_0$ to~$\s$, and let~$f\subset\s$ be a face.  Further
let~$H_{f,\g}$ be the subgroup of~$\M_r$ which is induced via~$m_r$ by
all elements of~$\pi_1(\Sph^d\setminus F, y_0)$ of the
form~$[\varphi^{-1}(\realPath{\g\,\delta\,\g^{-1}})]$, where~$\delta$
is any closed facet path in $\Star_{S_i}(f)$ based at~$\s$. The
subgroup~$H_{f,\g}$ has at least~$\Dim(f)+1$ trivial orbits, namely,
the orbits corresponding to the vertices of~$f$, and for~$g\subset f$
we have that the set of trivial orbits of~$H_{f,\g}$ contains the
trivial orbits of~$H_{g,\g}$. We consider the following three case:
\begin{enumerate}[(i)]
\item The intersection $\s\cap D_i$ is a ridge~$f$. Let~$\gamma$ be a facet path
  in~$D_i\cup\s$ from~$\s_0$ to~$\s$, and color~$\s$ (and
  hence~$f$) by the coloring induced along~$\g$ by the fixed
  coloring of~$\s_0$. Now keep the coloring of~$f$, but
  change the color of the remaining vertex~$v=\s\setminus f$ to any trivial orbit
  of~$H_{v,\g}$; see Figure~\ref{fig:case_i_ii} (right).

\begin{figure}[t]
  \centering
  \begin{overpic}[height=5.8cm]{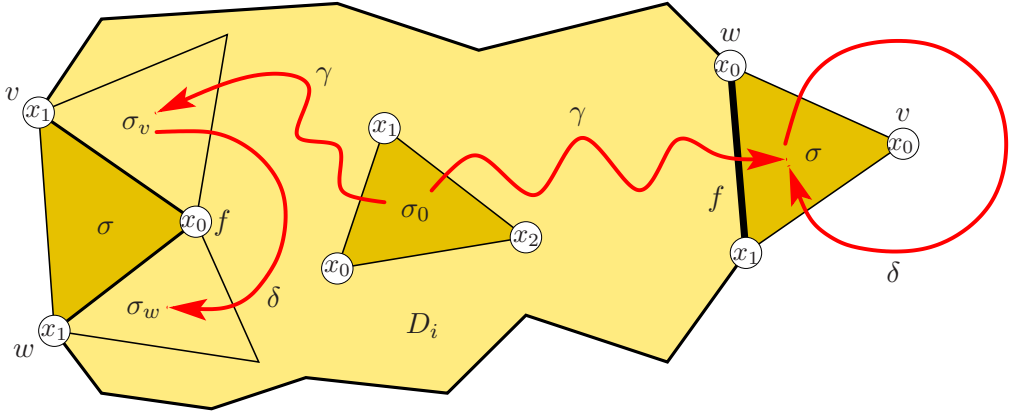}
    \put(38.5,20.5){$\s_0$}
    \put(31,15){\small$x_0$}
    \put(35.7,28.8){\small$x_1$}
    \put(49.5,18){\small$x_2$}
    \put(39,9){$D_i$}
    \put(55,30){$\g$}
    \put(86,14){$\delta$}
    \put(69.2,34.8){\small$x_0$}
    \put(86.1,27.1){\small$x_0$}
    \put(71,16.2){\small$x_1$}
    \put(78,26){$\s$}
    \put(68.5,21.5){$f$}
    \put(86.9,29.9){$v$}
    \put(69.7,37.8){$w$}
    \put(30.2,34){$\g$}
    \put(25.4,11.6){$\delta$}
    \put(17,19.5){\small$x_0$}
    \put(1.8,30.3){\small$x_1$}
    \put(3.3,8.7){\small$x_1$}
    \put(8.8,18.8){$\s$}
    \put(11.2,29.1){$\s_v$}
    \put(11.7,10.9){$\s_w$}
    \put(20.2,19.2){$f$}
    \put(-.2,32){$v$}
    \put(0.6,6.6){$w$}
  \end{overpic}
  \caption{Case~(i): The $2$-ball~$D_i$ with the facet~$\s_0$ colored via~$\imath$ by
    the preimage~$\{x_0,x_1,x_2\}$ of~$y_0$ and induced coloring of
    the ridge~$f$ on the right hand side of the figure.
    The vertex~$v$ is colored~$x_0$ if any element of~$\M_r$
    corresponding via~$m_r\circ\varphi^{-1}$ to a facet path of the
    form~$\g\,\delta\,\g^{-1}$ maps~$x_0$ to itself.
    Case~(ii): The induced coloring of the co-dimension
    2-face~$f$ and the vertices~$v$ and~$w$ on the left. The
    edge~$\{v,w\}$ is subdivided
    if the facet path~$\g\,\delta\,(\s_w,\s,\s_v)\,\g^{-1}$ corresponds
    via~$m_r\circ\varphi^{-1}$ to the identity in~$\M_r$.
    \label{fig:case_i_ii}}
\end{figure}

\item The intersection $\s\cap D_i$ equals two ridges~$f\cup v$
  and~$f\cup w$ with a common co-dimension 2-face~$f$.
  Let~$\s_v\in D_i$ be the facet intersecting~$\s$ in~$f\cup v$,
  let~$\s_w\in D_i$ be the facet intersecting~$\s$ in~$f\cup w$,
  and choose facet paths~$\gamma$ from~$\s_0$ to~$\s_v$ in~$D_i$
  and~$\delta$ from~$\s_v$ to~$\s_w$ in~$\Star_{D_i}(f)$.
  The fixed coloring of~$\s_0$ induces along~$\gamma$,
  respectively~$\g\,\delta$, colorings on~$f\cup v$
  and~$f\cup w$, and the colorings coincide on~$f$. Now we change
  the color of~$w$ according to
  $\m_r([\varphi^{-1}(\realPath{\g\,\delta\,(\s_w,\s,\s_v)\,\g^{-1}})])$,
  which is either a transposition (changing the color of~$w$)
  or the identity; see Figure~\ref{fig:case_i_ii} (left).
  
\item Otherwise set $S_{i+1}=S_i$ and let
  $(\s_0,\s_1,\dots,\s_{l_i},\s)$ be the marked beginning of a
  shelling of~$S_{i+1}$.
\end{enumerate}

We obtained a (possibly inconsistent) coloring of the vertices of~$\s$
in the cases~(i) and~(ii). Note that the coloring of~$\s$ induces a
consistent coloring on~$D_i\cap\s$, and that there is at most one
\emph{conflicting edge}~$\{v,w\}$, that is,~$v$ and~$w$ are colored
the same. A consistently colored
subdivision of~$\s$ is constructed in at most~$d-1$ subdivisions of~$\s$
with exactly one conflicting edge~$e$ each, where each subdivision is
obtained from the previous one by stellar subdividing~$e$:
Let~$f_e\subset\s$ be the unique minimal face such that $|e|\subset|f|$
holds and denote by~$C_e$ the set of trivial orbits of~$H_{f_e,\g}$. 
Now color the
new vertex~$v_e$ with an element of~$C_e$ which is not the color of any
vertex~$v_{e'}$ subdividing an edge~$e'$ with $f_{e'}\subset
f$. Note that~$C_e$ is the entire
preimage~$r^{-1}(y_0)$ if~$f_e$ is a co-dimension 1-face, and that~$C_e$
has at least one element distinct from the colors of all~$v_{e'}$
for~$f_{e'}\subset f_e$. If~$C_e$ contains
the one color~$x\in r^{-1}(y_0)$ not used in the coloring
of~$\s$, color~$v_e$ by~$x$ and terminate the subdivision process. 

\begin{figure}[htbp]
  \centering
  \begin{overpic}[height=5.7cm]{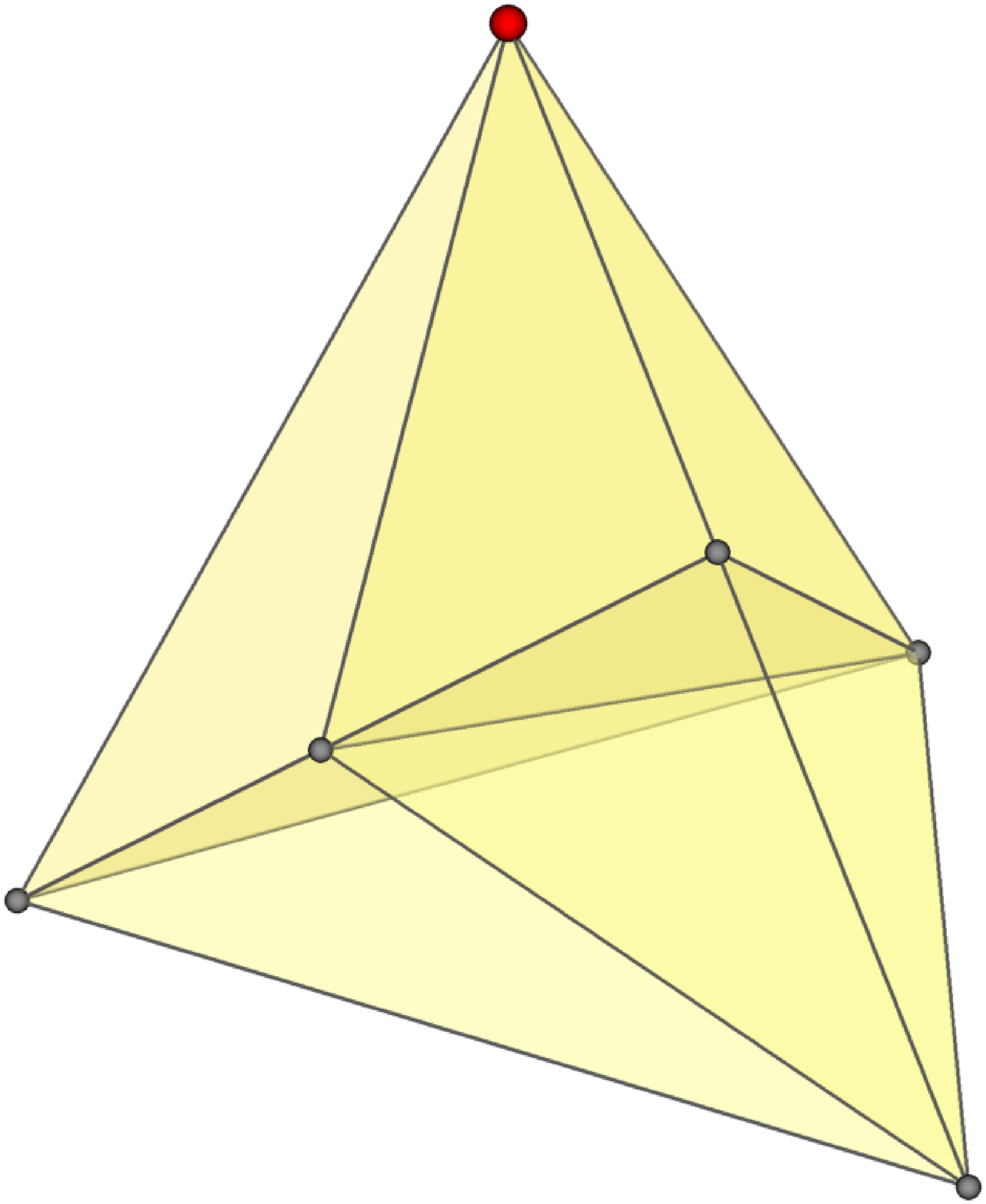}
    \put(22,96){$\{{\mathbf x_2}\}\,v$}
    \put(83,0){$\{{\mathbf x_2}\}$}
    \put(-12,23.5){$\{{\mathbf x_0}\}$}
    \put(79.1,43.9){$\{{\mathbf x_1}\}$}
    \put(61,54){\includegraphics[height=.25cm,angle=40]{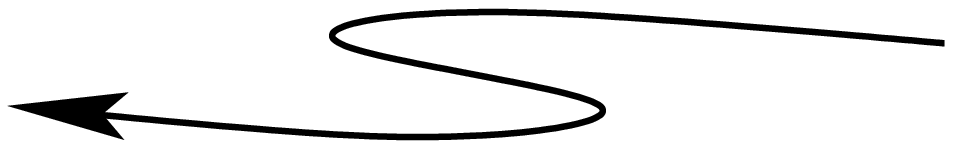}}
    \put(85.5,74){$\{{\mathbf x_0},x_2\}$}
    \put(1.5,40.1){\includegraphics[height=.25cm,angle=140]{arrow2}}
    \put(-18.3,59){$\{x_0,x_1,x_2,{\mathbf x_3}\}$}
  \end{overpic}
  \hspace{2.8cm}
  \begin{overpic}[height=5.7cm]{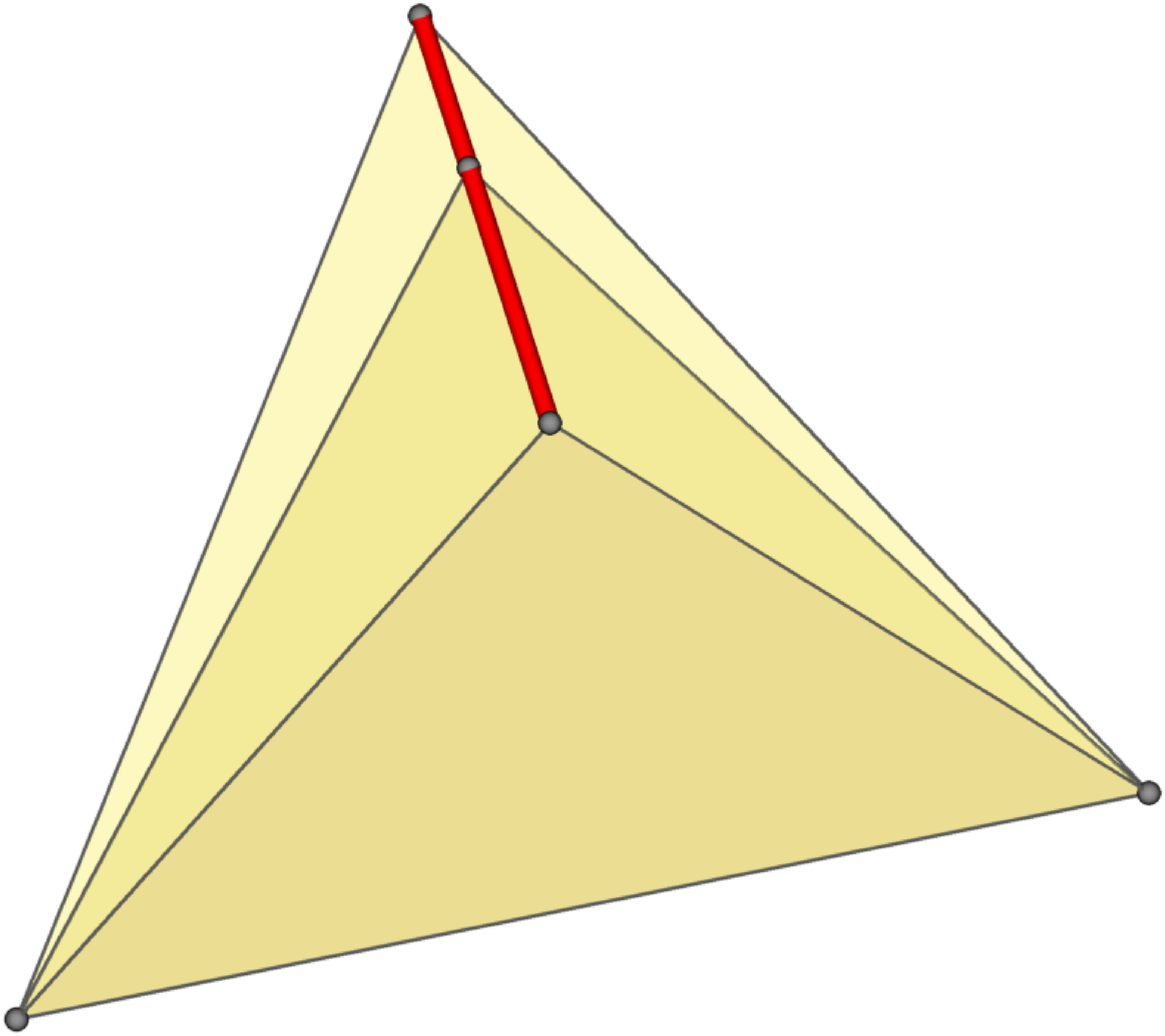}
    \put(17,87){$\{{\mathbf x_2}\}\,v$}
    \put(26,52.5){$\{{\mathbf x_2}\}\,w$}
    \put(-11,4.2){$\{{\mathbf x_0}\}$}
    \put(98.5,22.8){$\{{\mathbf x_1}\}$}
    \put(43.5,74.5){\includegraphics[height=.25cm]{arrow2}}
    \put(72.5,75.9){$\{x_0,x_1,x_2,{\mathbf x_3}\}$}
  \end{overpic}
  \caption{Coloring of the vertices of the refinement of~$\s$ in
    case~(i) (on the left) and case~(ii) (on the right). The minimal
    free face~$v$, respectively~$\{v,w\}$, is marked.  Each vertex~$v_e$
    is labeled by the trivial orbits of~$H_{f_e,\g}$ and the vertex color is
    printed bold.
    \label{fig:labeling}}
\end{figure}

This completes the construction of~$S_{i+1}$ in the cases~(i) and~(ii), and we define the
marked beginning of a shelling of~$S_{i+1}$ by
$(\s_0,\s_1,\dots,\s_{l_i})$ followed by the facets of the refinement
of~$\s$ in an appropriate order.

It remains to prove that the algorithm described above terminates
and that $p:\widehat{S}\to S$ is a branched cover equivalent
to~$r:X\to\Sph^d$.
Since~$S$ is shellable and hence nice,~$p$ is a branched cover by see
Corollary~\ref{cor:unf_are_covers}. The following
Lemmas~\ref{lem:equivalence} and~\ref{lem:non_free_edge} prove the
equivalence of~$p$ and~$r$, while termination of the construction
above is provided by Lemma~\ref{lem:terminates}.

\begin{lem}\label{lem:equivalence}
  The branched covers $p:\widehat{S}\to S$ and $r:X\to\Sph^d$ are equivalent.
\end{lem}

\begin{proof}
  In order to show the equivalence of the branched covers~$p$ and~$r$
  we prove by induction that the following holds for $0\leq i\leq l$:
  \begin{enumerate}[(I)]
  \item For any closed facet path~$\g$ based at~$\s_0$ in~$D_i$ we
    have
    \[
    \langle\g\rangle=\imath_*\circ\m_r([\varphi^{-1}(\realPath{\g})]).
    \]
  \item Let~$v\in D_i$ be a vertex, and let~$\gamma$ be a facet path in~$D_i$
    from~$\s_0$ to a facet~$\s$ containing~$v$.
    Then the color induced on~$v$ along~$\g$ by the fixed
    coloring of~$\s_0$ is a trivial orbit of~$H_{v,\g}$.
  \end{enumerate}
  
  We remark 
  that~(I) implies that~$\varphi$ restricted
  to~$\varphi^{-1}(|D_i|)$ is a homomorphism of pairs 
  $(\varphi^{-1}(|D_i|),F\cap\varphi^{-1}(|D_i|))\to(|D_i|,|\odd(D_i)|)$
  and that the Diagram~\eqref{eq:constructing_branched_covers} commutes.
  Finally,~(I) and~(II) are met for the pair~$(S_0,D_0)=(S',\s_0)$, and commutativity of
  Diagram~\eqref{eq:constructing_branched_covers} proves the
  equivalence of $r:X\to\Sph^d$ and $p:\widehat{S}\to S$
  for~$i=l$; see Theorem~\ref{thm:partial_unf}.

  We show that~(I) and~(II) hold for the pair~$(S_{i+1},D_{i+1})$
  provided they hold for the pair~$(S_i,D_i)$.  Recall that we denote
  the first facet~$\s_{l_i+1}$ of the shelling of~$S_i$ not contained
  in~$D_i$ by~$\s$. The simplicial complex~$D_i$ is contractible
  and hence~$\Pi_0(D_i,\s_0)=\Pi(D_i,\s_0)$ is generated by closed facet
  paths around (odd) co-dimension 2-faces by Theorem~\ref{thm:IJ_1}.
  Thus it suffices to verify~(I) for closed facet paths around
  (interior) co-dimension 2-faces by examining the three
  cases~(i),~(ii), and~(iii).
  \begin{enumerate}[(i)]
  \item The intersection $\s\cap D_i$ is a ridge~$f$. New interior
    co-dimension 2-faces in~$D_{i+1}$ arise only in the refinement
    of~$\s$, which is foldable by construction. Since~$\varphi(F)$ does not intersect
    the interior of~$|\s|$, any facet path around a new interior co-dimension 2-face
    corresponds to the identity of~$\M_r$ and~(I) holds by
    Lemma~\ref{lem:foldable_and_odd}.
  \item The intersection $\s\cap D_i$ equals two ridges~$f\cup v$
    and~$f\cup w$ with a common co-dimension 2-face~$f$.
    By induction hypothesis~(II) holds for the vertices of~$f$ in~$D_i$
    and thus~(I) follows for the new interior co-dimension
    2-face~$f$ of~$D_{i+1}$ by construction. As for any new interior
    co-dimension 2-face in the refinement
    of~$\s$,~(I) holds (as in case~(i)) since the refinement is
    foldable and~$\varphi(F)$ does not intersect
    the interior of~$|\s|$.
  \item Otherwise there is no co-dimension 2-faces~$f\subset\s$
    with a free corresponding edge~$e_f=\s\setminus f$ and~(I) follows from
    Lemma~\ref{lem:non_free_edge}.
  \end{enumerate}

  Having established~(I), it suffices to
  verify~(II) for a single facet path~$\g$ in~$D_{i+1}$ from~$\s_0$ to any facet
  containing a given vertex~$v$. Thus~(II) holds by choice of color for any vertex
  added to~$D_i$ in the construction of the pair~$(S_{i+1},D_{i+1})$.
\end{proof}

\begin{lem}\label{lem:non_free_edge}
  If~$f\in\s$ is a co-dimension $2$-face with a non-free
  corresponding edge $e_f=\s\setminus f$, then~{\rm(I)}
  holds for any closed facet path based at~$\s_0$ around~$f$ in~$D_{i+1}$.
\end{lem}

\begin{proof}
  Let~$\g\,\delta\,\g^{-1}$ be a closed facet path based at~$\s_0$
  around~$f$ in~$D_{i+1}$, where~$\delta$ is a closed path around~$f$
  in~$\Star_{D_{i+1}}(f)$.  Since~$\{v,w\}=e_f$ is a non-free edge,
  there is a facet path~$\delta'$ in~$D_i$ with~$\realPath{\delta'}$
  homotopy equivalent to $|\{f_e,f\cup v,f\cup w\}|$
  in~$|D_i|\setminus|\odd(D_i)|$, and we assume~$\delta$ and~$\delta'$
  to have the same orientation; see Figure~\ref{fig:non_free}. Note
  that the complex~$\{f_e,f\cup v,f\cup w\}$ itself is homotopy
  equivalent to~$\Sph^1$.

\begin{figure}[t]
  \centering
  \begin{overpic}[height=5.8cm]{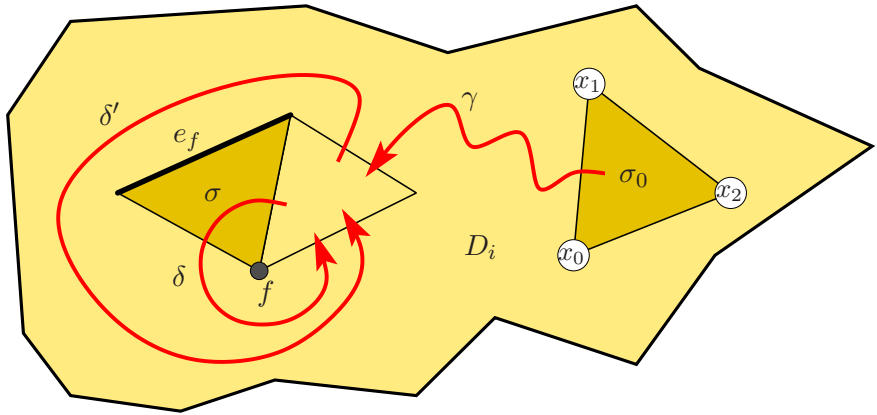}
    \put(63.8,18.4){\small$x_0$}
    \put(65.6,37.4){\small$x_1$}
    \put(81.1,25.4){\small$x_2$}
    \put(53.5,19){$D_i$}
    \put(70.5,27){$\s_0$}
    \put(25,25){$\s$}
    \put(30.8,13.9){$f$}
    \put(21.6,31.5){$e_f$}
    \put(53.2,35.5){$\g$}
    \put(21.5,15.5){$\delta$}
    \put(13.5,33.4){$\delta'$}
  \end{overpic}
  \caption{
    Case (iii): The
    paths~$\g$,~$\delta$, and~$\delta'$ if the corresponding
    edge~$e_f$ of a co-dimension 2-face~$f$ is non-free.
    \label{fig:non_free}}
\end{figure}

W.l.o.g. let~$m_r([\varphi^{-1}(\realPath{\g\,\delta\,\g^{-1}})])$ either be the identity or the
transposition~$(x_0,x_1)\in \M_r$. Each
transposition~$(x_i,x_j)$, for $i\not=j$, appears at most once in the
(unique) reduced representation of the element $a=m_r([\varphi^{-1}(\realPath{A})])\in\M_r$
corresponding to the facet path~$A=\g\,\delta'\,\g^{-1}$,
since~$A$ is composed from facet paths around co-dimension
2-faces of~$\s$.  Let $b=m_r([\varphi^{-1}(\realPath{B})])\in\M_r$ denote the element corresponding to
the facet path~$B=\g\,\delta'\,\delta^{-1}\,\g^{-1}$,
then~$a=(x_0,x_1)\circ b$ holds if and only if~$(x_0,x_1)$ is in the
reduced representation of~$a$, and we have~$a=b$ otherwise. Since~(I)
holds for~$D_i$ and hence in particular for the facet path~$A$, and
with
\[
A=\g\,\delta'\,\g^{-1}=\g\,\delta'\,\delta^{-1}\,\g^{-1}\,\g\,\delta\,\g^{-1} =B\,\g\,\delta\,\g^{-1},
\]
we conclude that the projectivity along~$\g\,\delta\,\g^{-1}$ is the
identity on the vertices of~$\s_0$ if and only if~$\g\,\delta\,\g^{-1}$
corresponds via~$\m_r\circ\varphi^{-1}$ to the identity in~$\M_r$, and
exchanges exactly the vertices colored~$x_0$ and~$x_1$ otherwise.
\end{proof}

The following Lemma~\ref{lem:terminates} proves termination of the
construction of the shellable $d$-sphere~$S$ and completes the proof of
Theorem~\ref{thm:constructing_branched_covers}.

\begin{lem}\label{lem:terminates}
  The shellable $d$-sphere~$S$ is obtained by finitely many stellar
  subdivisions of edges.
\end{lem}

\begin{proof}
  We prove that no facet will be subdivided more than a finite number
  of times in the construction of~$S$. The facet~$\s_{l_i+i}$ is
  subdivided at most~$d-1$ times in the construction~$S_{i+1}$
  from~$S_i$, and no facet in~$D_i$ is
  subdivided. The refinement of~$\s_{l_i+1}$ is added to~$D_i$ to
  define~$D_{i+1}$ and no facet in the refinement will be subdivided
  any further.

  Problems may accrue since subdividing~$\s_{l_i+1}$ results in
  subdividing other facets (not in~$D_i$)
  intersecting~$\s_{l_i+1}$, and each facet of the refinement of an
  intersecting facet appears in the shelling, yet
  is not in~$D_{i+1}$. Thus a facet might get subdivided over and
  over again.

  For a face~$f\in S'$ let~$L_{f,i}\subset S_i$ denote the refinement
  of~$f$ in~$S_i$.  W.l.o.g. we may assume that the facets of the
  refinement~$L_{\s,i}$ of any facet $\s\in S'$ appear consecutively
  in the shelling order of~$S_i$.  Let~$\s\in S'$ be a fixed facet and
  let~$i_0$ be the number such that~$S_{i_0}$ is the $d$-sphere
  with~$\s_{l_{i_0}+1}$ is the facet of~$L_{\s,{i_0}}$ appearing first
  in the shelling order, that is,~$S_{i_0+1}$ is constructed by adding
  (a refinement) of the first facet of~$L_{\s,{i_0}}$ to~$D_{i_0}$.
  Thus we obtain an induced coloring of the boundary vertices
  of~$L_{\s,{i_0}}$ which is consistent on $D_{i_0}\cap L_{\s,{i_0}}$
  by construction.  Since~$\varphi(F)$ does not intersect the interior
  of~$|L_{\s,{i_0}}|$ and by Lemma~\ref{lem:foldable_and_odd}, it
  remains to prove that this coloring of $D_{i_0}\cap L_{\s,{i_0}}$
  extends to a foldable refinement of~$L_{\s,{i_0}}$ obtained via a
  finite series of stellar subdivisions.

  Observe that each facet of~$L_{\s,{i_0}}$ is the cone over a
  $(d-1)$-simplex in the boundary of~$L_{\s,{i_0}}$
  and that~$L_{\s,{i_0}}$ has no interior vertices: This is obviously true
  for~$L_{\s,0}=\s$. For~$1\leq i\leq {i_0}$ let~$\Cone(f)$ be a facet
  of~$L_{\s,i-1}$ with~$f$ is a boundary $(d-1)$-simplex. Now
  if~$\Cone(f)$ is subdivided via stellarly subdividing an edge~$e\in
  f$, both facets replacing~$\Cone(f)$ are cones over boundary
  $(d-1)$-simplices
  which in turn are obtained from~$f$ by replacing one vertex of~$e$
  by the new vertex subdividing~$e$.

  We strengthen the statement above and claim that each
  facet of~$L_{\s,{i_0}}$ is the cone over a
  $(d-1)$-simplex in~$D_{i_0}\cap L_{\s,{i_0}}$. To this end note the
  trivial fact
  that if~$e\in L_{g,i}$ is an edge of the subdivision of a boundary $k$-face
  $g\in\s$ and let $\{f_j\}_{1\leq j\leq d-k}$ be the boundary $(d-1)$-faces of~$\s$
  with $g=\bigcap_{1\leq j\leq d-k} f_j$, then there is a $(d-1)$-simplex
  in each~$L_{f_j,i}$ containing~$e$. Thus if for
  some~$i<{i_0}$ an
  edge~$e$ is subdivided when adding the simplex~$\s_{l_i+1}$ to~$D_i$ which
  intersects~$L_{\s,i}$ in a low dimensional face, then at least
  one of the boundary $(d-1)$-simplices of~$L_{\s,i}$ containing~$e$
  will be added to~$D_{i'}\cap L_{\s,{i'}}$ at some point~$i<i'\leq i_0$.

  Returning to the consistent coloring of~$D_{i_0}\cap L_{\s,{i_0}}$
  we conclude that all vertices of~$L_{\s,{i_0}}$ are colored since
  there are no interior vertices, and that each facet~$\Cone(f)$
  of~$L_{\s,{i_0}}$ has at most one conflicting edge since the
  boundary $(d-1)$-simplex~$f\subset D_{i_0}\cap L_{\s,{i_0}}$ is consistently colored.
  Hence~$\Star_{L_{\s,{i_0}}}(e)$ of an conflicting edge~$e$ does not
  contain any other conflicting edges and we
  consider~$\Star_{L_{\s,{i_0}}}(e)$ independently.

  Now~$\Star_{L_{\s,{i_0}}}(e)$ is subdivided only finitely many
  times since~$H_{v,\g}$ is trivial for any new vertex~$v$ (except for finitely
  many vertices in the boundary of~$|\Star_{L_{\s,{i_0}}}(e)|$) and hence
  the construction (case~(i) and~(ii)) induces a linear order on
  the colors used to color the new vertices.
\end{proof}

\begin{rem}
  It appears as if the shellable $d$-sphere~$S$ may be constructed
  along a spanning tree of the dual graph~$\Gamma^*(S')$ instead of a
  shelling, though the construction would become substantially more
  complicated. Using a spanning tree of~$\Gamma^*(S')$ would eliminate
  the some how (to the theory of branched covers) alien concept of a
  shelling, and would allow for more general base spaces, e.g. PL
  $d$-manifolds.
\end{rem}

\begin{figure}[htbp]
  \centering
  \includegraphics[width=.48\textwidth]{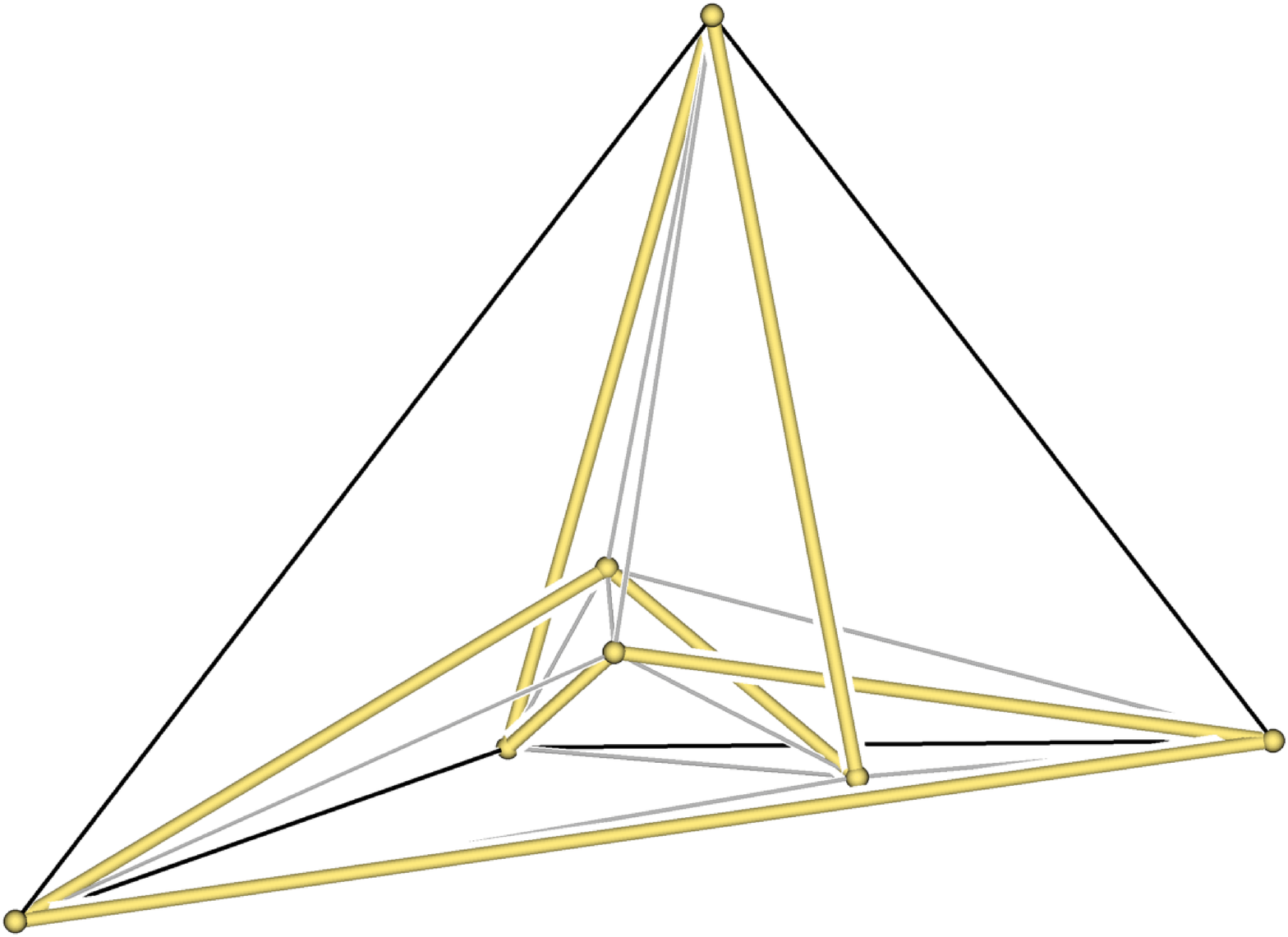}
  \hspace{.4cm}
  \includegraphics[width=.48\textwidth]{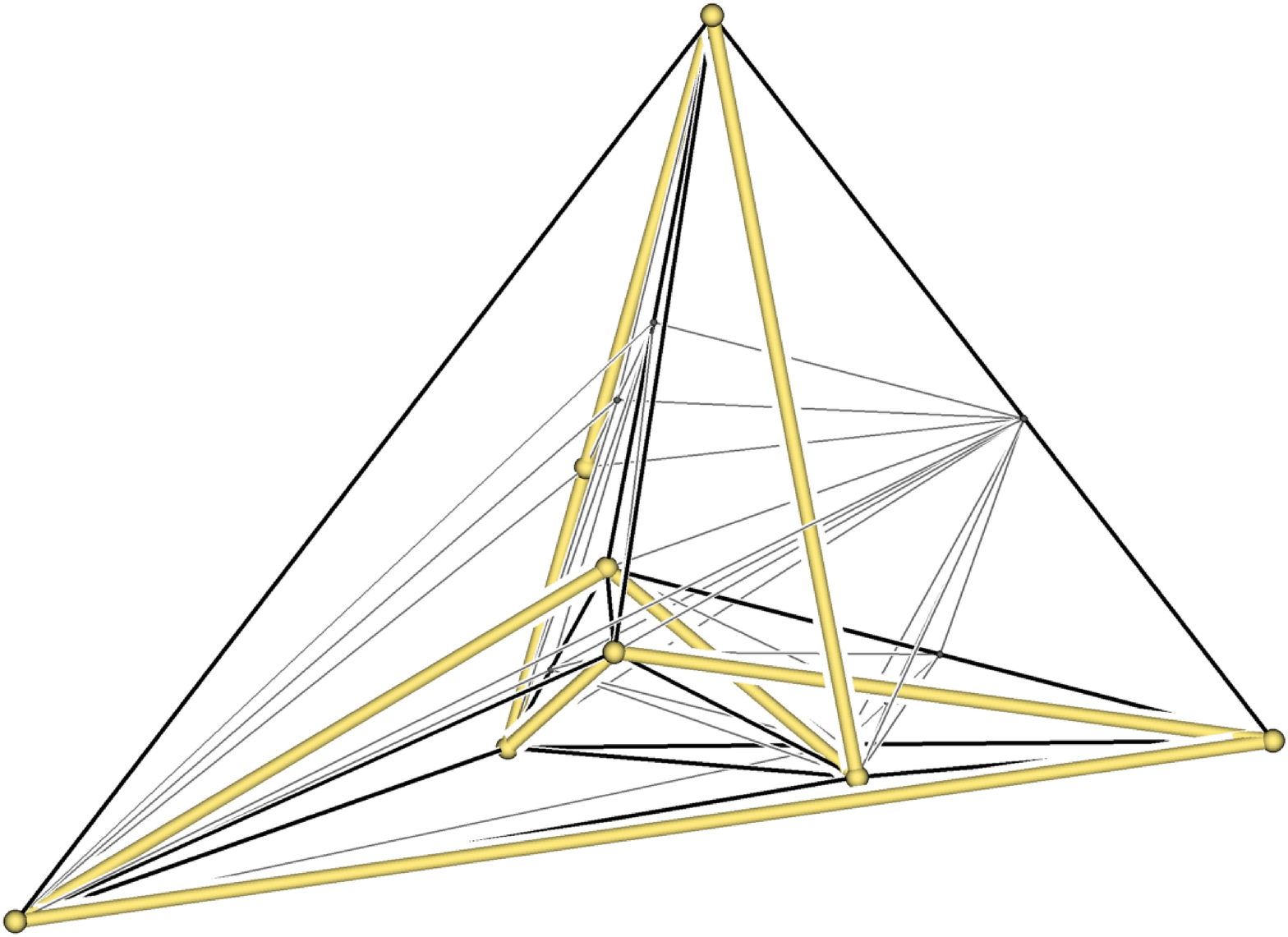}
  \caption{
    Constructing the trefoil as odd subcomplex of a
    $3$-sphere~$S$ with~$\Pi(S)$ isomorphic to the symmetric group on
    three elements. On the left the Schlegel diagram of the cyclic
    $4$-polytope~$C_{4,7}$ on seven vertices with the trefoil embedded
    in the $1$-skeleton. On the left~$S$ as a subdivision of the Schlegel diagram after stellarly
    subdividing eight edges of~$C_{4,7}$. The odd subcomplex is marked
    and the $f$-vector of~$S$ reads~$(15,63,96,48)$; \emph{watch}~\cite{witte:CYC}.
    \label{fig:trfoil_in_cyclic} }
\end{figure}

Applying Theorem~\ref{thm:constructing_branched_covers} to the results
of Hilden~\cite{hilden:TFBC} and Montesinos~\cite{montesinos:MFBC},
Piergallini~\cite{piergallini:FMS4}, and Iori~\&
Piergallini~\cite{iori_piergallini:4MFNS} we obtain the following
three corollaries.

\begin{cor}\label{cor:hilden_montesinos}
  Let $d=2$ or $d=3$. For every closed oriented $d$-manifold~$M$ there is a polytopal
  $d$-sphere~$S$ such that one of the connected
  components~$\widehat{S}$ of the partial unfolding
  of~$S$ is a combinatorial $d$-manifold homeomorphic to~$M$. The
  projection $\widehat{S}\to S$ is a simple $d$-fold branched cover
  branched over finitely many points for~$d=2$, respectively a link for~$d=3$.
\end{cor}

\begin{cor}\label{cor:piergallini}
  For every closed oriented PL $4$-manifold~$M$ there is a polytopal
  $4$-sphere~$S$ such that one of the connected
  components~$\widehat{S}$ of the partial unfolding
  of~$S$ is a combinatorial $4$-manifold PL-homeomorphic to~$M$. The
  projection $\widehat{S}\to S$ is a simple $4$-fold branched cover
  branched over a PL-surface with a finite number of cusp and node
  singularities.
\end{cor}

\begin{cor}\label{cor:iori_piergallini}
  For every closed oriented PL $4$-manifold~$M$ there is a polytopal
  $4$-sphere~$S$ such that the partial unfolding~$\widehat{S}$ of~$S$
  is a combinatorial $4$-manifold PL-homeomorphic to~$M$. The
  projection $\widehat{S}\to S$ is a simple $5$-fold branched cover
  branched over a locally flat PL-surface.
\end{cor}

A weaker version of Corollary~\ref{cor:hilden_montesinos} was already
established by Izmestiev~\& Joswig~\cite{izmestiev_joswig:BC} and
later by Hilden, Montesinos-Amilibia, Tejada~\& Toro~\cite{MR2218370}.
A weaker version of Corollary~\ref{cor:piergallini} can be found
in~\cite{witte:DISS}.

Stellar subdivision of an edge~$e\in S$ of a combinatorial $d$-manifold~$S$
changes the parity of the co-dimension 2-faces in~$\Link_S(e)$. Since
the link of an edge of~$S$ is a (combinatorial) $(d-2)$-sphere, we
obtain the following Corollary~\ref{cor:branching_set}. A
topological proof for arbitrary simple branched covers is avaliable
by Izmestiev~\cite{izmestiev:ADD}.

\begin{cor}\label{cor:branching_set}
  The branching set of a branched cover $r:X\to\Sph^d$ as described
  in Theorem~\ref{thm:constructing_branched_covers} is the symmetric difference of
  finitely many $(d-2)$-spheres.
\end{cor}

We conclude this section by a remark and a conjecture as to which
branched covers $r:X\to\Sph^d$ may be obtained via the method
presented above. In other words, which branching sets can be embedded
via a homomorphism $\varphi:\Sph^d\to|S'|$ into the co-dimension
$2$-skeleton of a shellable simplicial $d$-sphere~$S'$.

\begin{rem}
  For $d\geq6$ there are branching sets non-embedable into the co-dimension
$2$-skeleton of a shellable simplicial $d$-sphere: Freedman~\&
Quinn~\cite{freedman_quinn:TO4MF} constructed a 4-manifold which does
not have a triangulation as a combinatorial manifold. In fact, there
are 4-manifolds which can not be triangulated at all~\cite[p.~9]{lutz:MFV}.
\end{rem}

The branching set of a branched cover $r:X\to\Sph^d$ for
$d\leq5$ is at most 3-dimensional and since there is no difference in between PL
and non-PL topology up to dimension three, we conjecture the following.

\begin{conj}
  For $d\leq5$ every branched cover $r:X\to\Sph^d$ can be obtained via
  the partial unfolding of some polytopal $d$-sphere.
\end{conj}


\section{Extending Triangulations}
\label{sec:extending_triangulations}

\noindent
A first assault on how to extend triangulation and coloring is by Goodman~\&
Onishi~\cite{goodman_onishi:ETC}, who proved that a 4-colorable
triangulation of the 2-sphere may be extended to a 4-colorable
triangulation of the 3-ball. Their result was improved independently by
Izmestiev~\cite{izmestiev:EOC} and~\cite{witte:DIPL} to arbitrary
dimensions. Here we generalize the construction to arbitrary
simplicial complexes with $k$-colored subcomplexes.

\begin{thm}\label{thm:extending_triangulations}
  Given a simplicial $d$-complex~$K$ and a $k$-colored induced
  subcomplex~$L$, then there is a finite series of stellar subdivisions
  of edges, such that the resulting simplicial complex~$K'$ has a
  $\max\{k,d+1\}$-coloring,~$K'$ contains~$L$ as an induced
  subcomplex, and the $\max\{k,d+1\}$-coloring of~$K'$ induces the
  original $k$-coloring on~$L$.
\end{thm}

\begin{proof}
  We may assume~$K$ to be pure.  Let $K_0=K$ and assign~$0$ to all
  vertices not in~$L$. For $1\leq i\leq d$ we obtain the simplicial
  complex~$K_i$ from~$K_{i-1}$ by stellar subdividing all conflicting edges with
  both vertices colored~$i-1$ in an arbitrary order. The new vertices
  are colored~$i$. We prove by induction that for $0\leq j\leq i-1$
  and each facet~$\s\in K_i$ there is exactly one vertex $v_j\in\s$
  colored~$j$. The assumption holds for~$K_0$ and completes
  the proof for $K'=K_d$. Note that since~$L$ is properly colored, no
  edges in~$L$ are subdivided and~$L$ is an induced subcomplex of
  any~$K_i$ for~$0\leq i\leq d$.

  To prove the induction hypothesis for~$K_i$, we again use an inductive
  argument: Let~$\s$ be a facet of a subdivision of~$K_{i-1}$ produced
  in the making of~$K_i$. Assume that each color less than~$i-1$ appears
  exactly once in~$\s$, and let~$l\geq2$ be the number of
  $(i-1)$-colored vertices of~$\s$. This assumption clearly holds for
  any facet of~$K_{i-1}$ for some $l\leq d-i+2$. After subdividing an
  $(i-1)$-colored conflicting edge of~$\s$ and assigning the color~$i$ to the
  new vertex, each of the two new facets has~$l-1$ vertices
  colored~$i-1$, and each color less than~$i-1$ appears exactly once. Thus any
  facet of~$K_{i-1}$ has to be subdivided into at most $2^{d-i+1}$ simplices in
  order for~$K_i$ to meet the induction hypothesis.
\end{proof}

Izmestiev gives a result similar to
Theorem~\ref{thm:extending_triangulations} in~\cite{izmestiev:EOC},
but the following Remark~\ref{rem:extending_triangulations} points out
the advantage of using only stellar subdivisions of edges.

\begin{rem}\label{rem:extending_triangulations}
  Since only stellar subdivisions of edges are used to construct~$K'$
  from~$K$ all properties invariant under these subdivisions are
  preserved, e.g.  polytopality, regularity, shellability, and others.
  In the case that~$L$ is not induced, stellarly subdivide all
  edges~$\{v,w\}\in K\setminus L$ with~$v,w\in L$.  In order to obtain
  a small triangulation, one can try to (greedily) $(d+1)$-color a
  (large) foldable subcomplex first.
\end{rem}

\begin{cor}\label{cor:odd_of_mf}
  The odd subcomplex of a closed combinatorial $d$-manifold is the symmetric
  difference of finitely many $(d-2)$-spheres.
\end{cor}

\begin{figure}[htbp]
  \centering
  \psfrag{0}{\small\;0}
  \psfrag{1}{\small\;1}
  \psfrag{2}{\small\;2}
  \includegraphics[height=5cm]{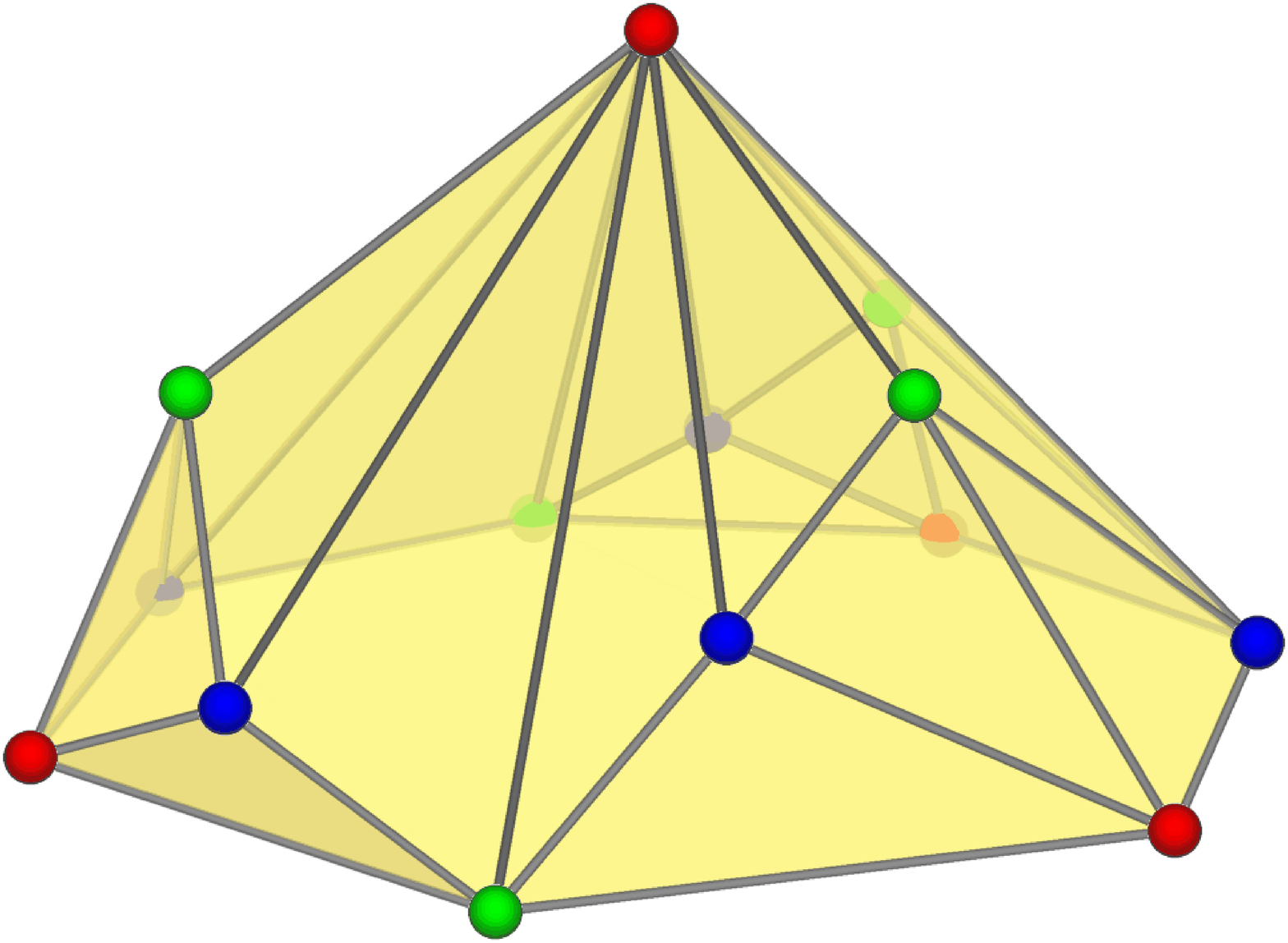}
  \hspace{1.2cm}
  \includegraphics[height=5cm]{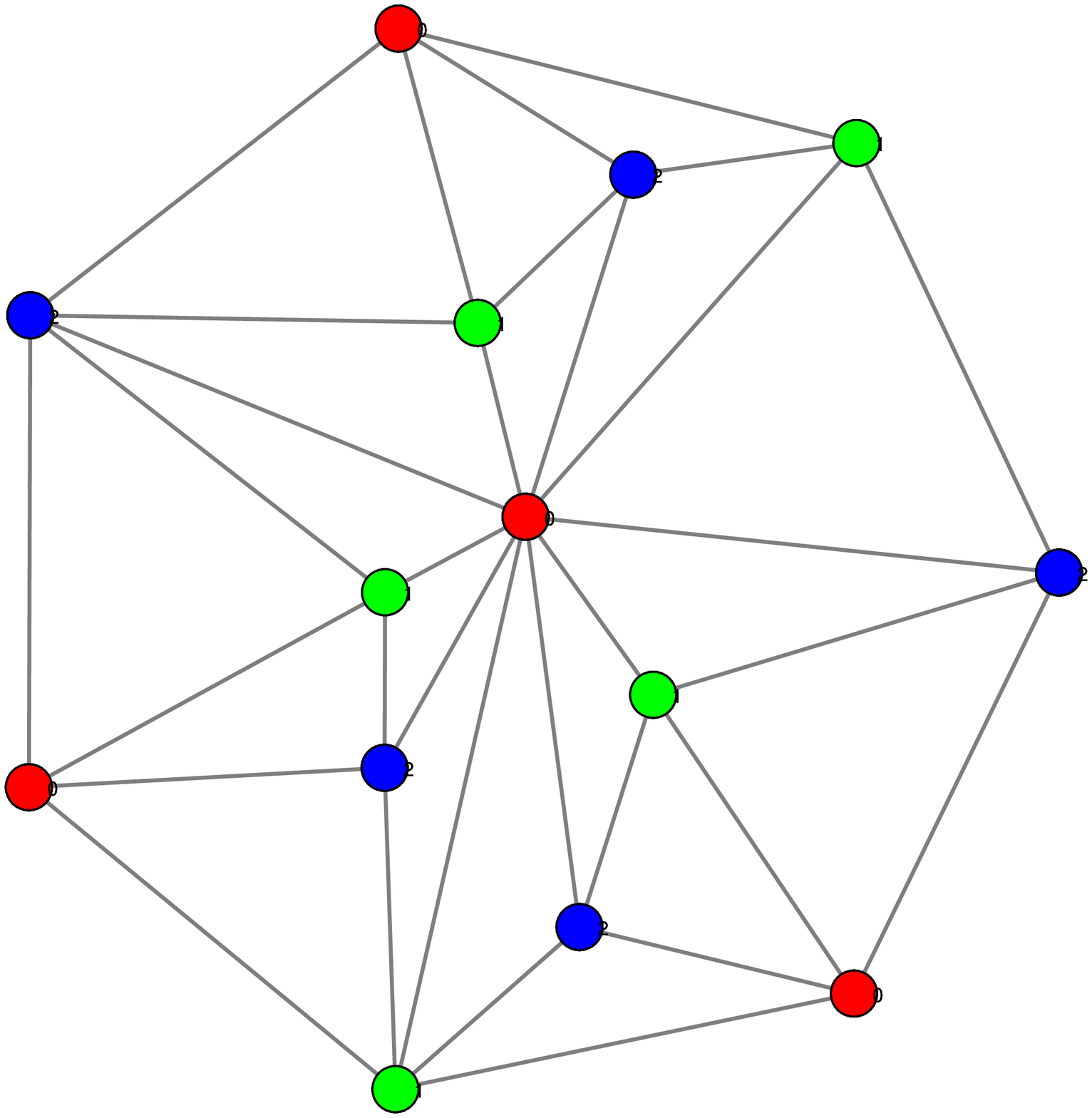}
  \caption{Convex hull of the extended triangulation of a 3-colored 7-gon
    and its Schlegel diagram.
    \label{fig:ext_coloring}}
\end{figure}

\begin{cor}\label{cor:extending_triangulations}
  Given a $k$-colored simplicial $(d-1)$-sphere~$S$, then there is a
  simplicial $d$-ball~$D$ with boundary equal to~$S$, such that there
  is a $\max\{k,d+1\}$-coloring of~$D$ which induces the original
  $k$-coloring on~$S$.
  The $d$-ball~$D$ is obtained from $\Cone(S)$ by a finite series of
  stellar subdivision of edges. In particular~$D$ is
  a combinatorial $d$-ball if~$S$ is a combinatorial $(d-1)$-sphere,
  shellable if~$S$ is shellable, and
  regular if~$S$ is polytopal; see Figure~\ref{fig:ext_coloring}.
\end{cor}

\begin{rem}
  Similar results as Corollary~\ref{cor:extending_triangulations} may
  easily be obtained for partial triangulations of CW-complex and
  relative handle-body decomposition of PL-manifold (with boundary).
\end{rem}

The partial unfoldings of two homeomorphic simplicial complexes~$K$
and~$K'$ need not to be homeomorphic in general. We present a notion
of equivalence of simplicial complexes which agrees with their
unfolding behavior, that is, we give sufficient criteria such that if
$p:\widehat{K}\to K$ and $p':\widehat{K'}\to K'$ are branched covers,
then~$p$ and~$p'$ are equivalent.

Assume~$K$ and~$K'$ to be strongly connected and that the odd
subcomplexes~$K_{\odd}$ and~$K'_{\odd}$ are \emph{equivalent}, that
is, there is a homeomorphism of pairs $\varphi:(|K|,|K_{\odd}|)\to
(|K'|,|K'_{\odd}|)$. Let $\s_0\in K$ be a facet, and~$y_0$ the
barycenter of~$\s_0$, and assume that the image $y'_0=\varphi(y_0)$ is
the barycenter of~$|\s'_0|$ for some facet $\s'_0\in K'$. Now~$K$
and~$K'$ are \emph{color equivalent} if there is a bijection
$\psi:V(\s_0)\to V(\s'_0)$, such that
\begin{equation}\label{equ:color_equivalent}
  \psi_*\circ\h_K=\h_{K'}\circ\varphi_*
\end{equation}
holds, where the maps
$\varphi_*:\pi_1(|K|\setminus|K_{\odd}|,y_0)\to\pi_1(|K'|\setminus|K'_{\odd}|,y'_0)$
and $\psi_*: \Sym(V(\s_0))\to \Sym(V(\s'_0))$ are the group
isomorphisms induced by~$\varphi$ and~$\psi$, respectively.

Observe that this is indeed an equivalence relation. The name ``color
equivalent'' suggests that the pairs $(K,K_{\odd})$ and
$(K',K'_{\odd})$ are equivalent, and that the ``colorings'' of~$K_{\odd}$ and~$K'_{\odd}$
by the $\Pi(K)$-action, respectively $\Pi(K')$-action, of projectivities
around odd faces are
equivalent. Lemma~\ref{lem:color_equiv} justifies this name.

\begin{lem}\label{lem:color_equiv}
  Let~$K$ and~$K'$ be color equivalent nice simplicial complexes. Then the
  branched covers $p:\widehat{K}\to K$ and $p':\widehat{K'}\to K'$ are
  equivalent.
\end{lem}

\begin{proof}
  With the notation of Equation~\eqref{equ:color_equivalent} we have that
  \[
  \xymatrix{
    & \pi_1(|K|\setminus|K_{\odd}|,y_0) \ar[r]^{\varphi_*}
    \ar[dl]_{\m_p} \ar[d]^{\h_K}&
    \pi_1(|K'|\setminus|K'_{\odd}|,y'_0)
    \ar[dr]^{\m_{p'}} \ar[d]^{\h_{K'}}& \\
    \M_p \ar[r]^{\imath_*} & 
    \Pi(K,\s_0)\ar[r]^{\psi_*} &
    \Pi(K',\s'_0)& 
    \M_{p'} \ar[l]_{\imath'_*}
  }
  \]
  commutes, since the Diagram~\eqref{equ:partial_unf} commutes and
  Equation~\eqref{equ:color_equivalent} holds. Theorem~\ref{thm:bc} completes the proof.
\end{proof}

\begin{prop}\label{prop:odd_by_stellar}
  For every strongly connected simplicial complex~$K$ there is a
  simplicial complex~$K'$ obtained from a foldable simplicial complex
  via a finite series of stellar subdivision of edges, such that~$K$
  and~$K'$ are color equivalent.
\end{prop}

Theorem~\ref{thm:constructing_branched_covers} proves
Proposition~\ref{prop:odd_by_stellar} above for shellable spheres.
We will not prove the general case and only give a sketch of the
construction for general~$K$.

Let~$L$ be a foldable simplicial complex obtained from~$K$ via a
finite series of stellar subdivisions according to
Theorem~\ref{thm:extending_triangulations}, that is, there is a
series $(K=K_0,K_1,\dots,K_l=L)$ where~$K_i$ is obtained
from~$K_{i-1}$ by stellar subdividing a single edge~$e_{i-1}\in K_{i-1}$. The idea is
to reverse the effect of the stellar subdivisions by subdividing
each edge~$e$ a second time in the reversed order, since stellar
subdividing~$e$ twice yields the anti-prismatic subdivision of~$e$
(which does not alter the color equivalence class).

We construct~$K'$ from~$L$ inductively in a series
$(L=L_l,L_{l-1},\dots,L_0=K')$ of simplicial complexes, where~$L_i$
is obtained from~$L_{i+1}$ by a finite series of stellar subdivision
of edges. The complexes~$L_i$ and~$K_i$ are color
equivalent: For a
facet path $(\s'_j)_{j\in J}$ in~$L_i$ associate the facet
path~$(\s_j)_{j\in J}$ in~$K_i$, where~$\s_j$ is the unique
facet, such that~$|\s'_j|$ lies in~$|\s_j|$.

We fix some notation in order to describe the construction of~$L_i$
from~$L_{i+1}$. Subdividing the edge~$e_i\in K_i$ in order to
construct~$K_{i+1}$ replaces~$e_i$ by two edges in~$K_{i+1}$, and we
call one of these two edges~$e'_i$. A facet
in~$\Star_{K_{i+1}}(e'_i)$ might get subdivided further in the
process of constructing
$K_{i+2},K_{i+3},\dots,K_l=L_l,L_{l-1},\dots,L_{i+1}$, and we
define~$L_{e'_i}$ as the subcomplex of~$L_{i+1}$ which
refines~$\Star_{K_{i+1}}(e'_i)$.

Note that~$e'_i$ is an edge of~$L_{e'_i}$, and that~$L_{e'_i}$
and~$\Star_{K_{i+1}}(e'_i)$ are color equivalent.  It follows that
the group of projectivities of~$L_{e'_i}$ has at least two trivial
orbits corresponding to the vertices of~$e'_i$.  Now~$L_i$ is
obtained from~$L_{i+1}$ by stellarly subdividing all edges with
vertices belonging to the same two trivial orbits as the vertices
of~$e'_i$.



\pagebreak
\bibliographystyle{mod_siam}
\bibliography{literature}

\end{document}